\documentclass[12pt,a4,oneside,reqno]{amsart} 
\usepackage{amsmath, amsthm}
\usepackage[psamsfonts]{amssymb}
\usepackage{color}
\usepackage[width=15cm,
            height=22cm]
           {geometry}
\usepackage{graphicx}

\usepackage[utf8]{inputenc}
\usepackage{fancyhdr}
\usepackage{xcolor}
\newcommand{\ab}{[a,b\hskip0.5mm]}
\newcommand{\eps}{\varepsilon}
\newcommand{\dd}{\hskip0.2mm\mbox{\rm d}}
\newcommand{\N}{{\mathbb N}}
\newcommand{\R}{{\mathbb R}}
\newcommand{\var}{\operatorname{var}}

\newcommand{\g}{[\hskip0.2mm\mbox{\rm d}g]}
\newcommand{\f}{[\hskip0.2mm\mbox{\rm d}f]}
\newcommand{\F}{[\hskip0.2mm\mbox{\rm d}F]}
\newcommand{\wt}{\widetilde}

\newcommand{\balfa}{{\boldsymbol{\alpha}}}
\newcommand{\ksi}{{\boldsymbol{\xi}}}

\newtheorem{theorem}{Theorem}[section]
\newtheorem{lemma}[theorem]{Lemma}
\newtheorem{definition}[theorem]{Definition}

\newtheorem{proposition}[theorem]{Proposition}
\newtheorem{corollary}[theorem]{Corollary}
\newtheorem{remark}[theorem]{Remark}

\def\skipaline{\removelastskip\vskip10pt plus 1pt minus 1pt}
\def\skiphalfaline{\removelastskip\vskip5pt plus 1pt minus 1pt}

\numberwithin{equation}{section}

\begin{document}
\title[]
{Role of the Harnack Extension Principle in the Kurzweil-Stieltjes Integral}

\makeatletter
\@namedef{subjclassname@2020}{%
  \textup{2020} Mathematics Subject Classification}
\makeatother


\author{Umi Mahnuna Hanung}
\address{Korteweg$-$de Vries Institute for Mathematics, Universiteit van Amsterdam, Amsterdam, Netherlands}
\email{H.UmiMahnuna@uva.nl}

\address{Department of Mathematics, Universitas Gadjah Mada, Sekip Utara Bulaksumur,
         55281 Yogyakarta, Indonesia}
\email{hanungum@ugm.ac.id}

\address{Institute of Mathematics, Czech Academy of Sciences,
          \v{Z}itn\'a 25, 11567 Praha 1, Czech Republic}

\thanks{The author was supported by the Ministry of Education, Culture, Research and Technology, Republic of Indonesia and and has been farther co-funded by the Grant No. 152/J01.1.28/PL.06.02/2022 of FMIPA UGM.}

\date{\today}

\begin{abstract}
Various kinds of Stieltjes integrals using gauge integration have become highly popular in the field of differential equations and other applications. In the theories of integration and of ordinary differential equations, convergence theorems provide one of the most widely used tools.
The Harnack extension principle, which discusses a sufficient condition for Kurzweil-Henstock integrable
functions on particular subsets of $(a,b)$ to be integrable on $\ab$ (see e.g., Theorem \ref{harnack}), is a key step to supply convergence theorems.
The Kurzweil-Stieltjes integral reduces to the Kurzweil-Henstock integral whenever
the integrator is an identity function. In general, if the integrator $F$ is
discontinuous on $[c,d]\subset\ab,$ then the values of the Kurzweil-Stieltjes integrals
\[
  \int_{c}^{d}\F\,g,\quad \int_{[c,d]}\F\,g,\quad
  \int_{[c,d)}\F\,g,\quad \int_{(c,d]}\F\,g,
  \mbox{ \ and \ } \int_{(c,d)}\F\,g
\]
need not coincide (see \cite[Section 5]{MHT}). Hence, the Harnack extension
principle in the Kurzweil-Henstock integral cannot be valid any longer for the Kurzweil-type Stieltjes integrals with discontinuous integrators.
The new concepts of equi-integrability and equiregulatedness involving elementary sets are pivotal to the notion of the Harnack extension principle for the Kurzweil-Stieltjes integration.

Moreover, in general, the existence of the integral
$\int_a^b\F\,g$ does not (even in the case of the identity integrator $F(x)=x$)
always imply the existence of the integral $\int_{T}\F\,g$ for every subset $T$
of $\ab.$ This follows from the well-known fact that, if e.g., $T\subset\ab$ is not
measurable, then the existence of the Lebesgue integral $\int_a^b g\,\dd t$
(which is a particular case of the Kurzweil-Henstock one) does not imply that the integral $\int_T g\,\dd t$ exists.
Therefore, besides having an interest in constructing the Harnack extension principle for the Kurzweil-Stieltjes integral,
the aim of this paper is also to demonstrate its role in guaranteeing the existence of the integrals $\int_{T}\F\,g$ for the arbitrary subsets $T$ of an elementary set $E.$
\end{abstract}

\maketitle

\section{Introduction}
One of the meaningful discussions in the topic of Kurzweil-Henstock integration concerns
the Harnack extension principle and Cauchy property (see e.g.,
\cite[Corollaries 7.10 and 7.11]{L1} and \cite[Theorems 1.4.6, 1.4.8 and 4.4.4]{L2}).
The Cauchy property was first used for the Riemann integral to integrate
functions unbounded in the neighborhood of a finite number of points (see e.g.,
\cite[Theorems 2.12 - 2.16]{R}). A similar idea has also been applied to integrate in
a Lebesgue sense functions not summable in the neighborhood of some points. Based
on the Cauchy property, C. G. A. Harnack suggested a method to calculate the integrals of functions
defined on an open set. The Cauchy property in the integral theory presents a sufficient
condition for the integrable functions on each $[c,d]\subset(a,b)$ to be integrable
on $\ab$ (see e.g., \cite{K1}, \cite{L1}, \cite{LPY2}, \cite{L2}). In the setting of the
Kurzweil-Henstock integral for real-valued functions, the Harnack extension principle
reads as follows (see e.g., \cite[Theorem 9.22]{G2}, \cite[Corollary 7.11]{L1},
\cite[Theorem 4.4.4]{L2}):\smallskip

\begin{theorem}\label{harnack}
Let $T\,{\subset}\,\ab$ be a closed set and let $\{[a_i,b_i]\,{:}\,i\in\N\}$
be a collection of pairwise disjoint intervals such that
$(a,b)\setminus T=\bigcup_{i=1}^{\infty}(a_i,b_i).$ Then, if $g$ is a real-valued function and the Kurzweil-Henstock integrals $\int_a^b g\chi_{T}\,\dd t$
and $\int_{a_i}^{b_i}g\,\dd t$ exist for all $i\in\N$ and the series
\begin{align*}\label{harnack2}
  \sum_{i=1}^{\infty}
     \sup\left\{\left|\int_r^t g\,\dd t\right|\,{:}\,a_i\le r\le t\le b_i\right\}
\end{align*}
converges, then the Kurzweil-Henstock integral $\int_a^b g\,\dd t$ exists
and
\[
  \int_a^b g\,\dd t
  =\int_a^b g\chi_{T}\,\dd t+\sum_{i=1}^{\infty}\int_{a_i}^{b_i} g\,\dd t.
\]
\end{theorem}\smallskip

The original definition of the gauge based integral, generally non-separated
functions of two variables, was
given by J. Kurzweil in the late 1950's and published
in his paper on generalized ordinary differential equations as an alternative definition
of the Perron/Denjoy integral (see e.g., \cite{K}, \cite{K4}). In the early 1960's,
R. Henstock independently rediscovered the analogous
definition of integral, and developed it into a systematic theory
(see e.g., \cite{RH1}, \cite{RH2}, \cite{RH3}).
Nowadays, the integral is known as the Kurzweil-Henstock
(or Henstock-Kurzweil) integral (see e.g., \cite{Fed_M_S}, \cite{KS}, \cite{K3}, \cite{SaPe_PB_V_JJO}).
Its definition based on Riemannian-type sums and refinements controlled by
gauges, which are also known as the generalized Riemann integral or gauge integral,
leads to a non-absolutely convergent integral that is more powerful than the Lebesgue integral, and also contains a special case, i.e., the Stieltjes-type integrals.
Throughout this paper, we work with the Kurzweil-Stieltjes integrals.
The simplest integral of this type is the Riemann-Stieltjes integral of the form
$\int_a^b\f\,g,$ in which a function $g{:}\,\ab\to\R$ called the integrand is
integrated with respect to another function $f{:}\,\ab\to\R$ referred to as the integrator
(see e.g., \cite{G2}, \cite{M3}, \cite{P}, \cite{St}). This integral
appeared for the first time in a famous treatise \cite{St} by T. J. Stieltjes.
Up to now, many authors have considered various kinds of Stieltjes integrals using the gauge integration
(see e.g., \cite{UMH1}, \cite{H2}, \cite{M1}, \cite{M3}, \cite{P}, \cite{S}, \cite{T1}), which have become highly popular in the field of differential equations and other applications (see e.g., \cite{Arr1}, \cite{Bon_MFed_Mesq}, \cite{Carl_Hei}, \cite{Fed_M_S},
\cite{FMS}, \cite{Krejci2}, \cite{Krejci3}, \cite{Leon_MOPE}, \cite{LiuW}, \cite{SaPe_MT}, \cite{SL}, \cite{YeG_ZLZ}).
In the literature, these integrals are known under several different names (e.g., Henstock-Stieltjes,
Perron-Stieltjes, and generalized Riemann-Stieltjes). In our opinion, all of these integrals
are special cases of the Kurzweil integral referred from \cite{K} or \cite{K4}.
Therefore, we prefer to call this integral the \emph{Kurzweil-Stieltjes
integral}.
\smallskip

The Kurzweil-Henstock integral has been generalized in various ways. For instance,
S. S. Cao (\cite{C}) noticed that the Kurzweil's definition can be easily extended to functions
with values in Banach spaces and investigated some of the properties of the abstract Kurzweil-Henstock
integral. This abstract Kurzweil-Henstock integral received further attention, such as
the monograph by \v{S}. Schwabik and G. Ye \cite{S5Y} that discusses these types of integrals, i.e.,
the McShane, Bochner, Dunford, and Pettis integrals for Banach space-valued functions, and compares the relationship between these various integrals.
Moreover, the fundamental results concerning
the Kurzweil-Stieltjes integral for Banach space-valued function integrals were given by \v{S}. Schwabik
in \cite{S} and \cite{S5}, where he called it the {\em abstract Perron-Stieltjes integral}.
The results obtained by \v{S}. Schwabik have undergone expansion by G. A.~Monteiro and M.~Tvrd\'y, completing them
to theory, such that it was applicable to prove some of the results on the continuous dependence of solutions
to generalized linear differential equations in a Banach space (see \cite{M1} and \cite{M1T}).
\smallskip

Convergence theorems for integrals concern the possibility of interchanging the limit and the integral
(see e.g., \cite{B1}, \cite{G3}, \cite{UMH1}, \cite{H1}, \cite{LPY2}, \cite{P}, \cite{S2}). The extension of the Kurzweil-Stieltjes
integral to the integration over elementary sets, i.e., sets that are finite unions of bounded intervals,
was presented in \cite[Section 5]{MHT}, where it was a useful ingredient for proving the bounded convergence
theorem for the abstract Kurzweil-Stieltjes integral. However, based on
\cite[Theorems 5.8, 5.9 and 5.10, Remark 5.12, Theorem 5.13]{MHT}, the Harnack extension
principle for the Kurzweil-Henstock integral, see e.g., Theorem \ref{harnack}, cannot be easily extended
to the Stieltjes-type integrals as, whenever the integrator $F$ is not continuous on $\ab,$
for a subinterval $J\subset\ab$ having an infimum and a supremum $c$ and $d,$ respectively,
the integrals
\begin{equation*}
  \int_{c}^{d}\F\,g,\quad \int_{[c,d]}\F\,g,\quad
  \int_{[c,d)}\F\,g,\quad \int_{(c,d]}\F\,g, \mbox{ \ and \ }
  \int_{(c,d)}\F\,g
\end{equation*}
need not have the same values, even if they all exist (see \cite[Remark 5.12]{MHT}). The reason is that the Stieltjes integrals over a one-point set need not be zero (see e.g., \cite[Proposition 5.7]{MHT}).
\smallskip

From the study of convergence theorems for gauge-type integrals,
the notion of
\textcolor[rgb]{0.00,0.00,0.00}{equi-integrability}
appeared,
whose idea is that there exists a single gauge $\delta$ that
works for all the functions in a sequence (see e.g., \cite{BoPi}, \cite{G3}, \cite{K5J}, \cite{M3}, \cite{M4}, \cite{S4V}, \cite{S5Y}).
Aside from extending the definition of the Kurzweil-Stieltjes integral
over arbitrary bounded sets (see Section \ref{IOABS}), to deal with the Harnack extension principle
for the Kurzweil-Stieltjes integral based on \cite{MHT}, it is necessary to develop the notion
of equi-integrability for Banach space-valued functions and investigate their fundamental properties,
including some important results regarding equiregulatedness, presented in Section \ref{equi}.
Furthermore, the theory in Sections \ref{IOABS} and \ref{equi} leads to a new Harnack
extension principle for the Kurzweil-Stieltjes integral, which significantly improves the results
from \cite{MHT}. Therefore, the goal of this paper is to provide sufficient conditions vouching the Harnack
extension principle for the Kurzweil-Stieltjes integral and then to show that it is applicable to integrate
functions over arbitrary closed subsets of an elementary set, as shown in Section \ref{HarnackApp}.

\skipaline\vskip-0.5mm

\section{Preliminaries}
In this section we will recall some terminologies and notations commonly used in
the literature.

Let $X,\,Y,$ and $Z$ be Banach spaces. The symbols
$\|\cdot\|_X,\,\|\cdot\|_Y,$ and $\|\cdot\|_Z$ stand for the norm in $X,\,Y,$ and $Z,$ respectively.

If there are bilinear mapping $B\,{:}\,X\times Y\to Z$ and $\beta\in[0,\infty)$ such that
\begin{align*}
  \big\|B(x,y)\big\|_Z\le\beta\,\|x\|_X\,\|y\|_Y \quad\mbox{for \ } x\in X,\, y\in Y,
\end{align*}
then the triple $(X,\,Y,\,Z)$ is a {\em bilinear triple} with respect to $B.$
In such a case, we write $\mathcal{B}=(X,Y,Z)$ and use the abbreviation $x\,y$ for $B(x,y).$
Besides a classical situation with $X=Y=Z=\R,$ a typical nontrivial example is, e.g.,
$\mathcal{B}=(\mathcal{L}(X,Z),X,Z),$ where $\mathcal{L}(X,Z)$ is the space of all linear
bounded operators $L\,{:}\,X\to Z,$ whereas $B(L,x)=L\,x\in Z$ for $x\in X$ and
$L\in\mathcal{L}(X,Z).$ Clearly, without any loss of generality, we may
assume that $\beta=1.$
\smallskip

Two intervals in $\R$ are said to be {\em disjoint} if their intersection is empty,
whereas they are said to be {\em non-overlapping} if their intersection contains at most one
point. In this study, with an {\em elementary set}, we understand a finite union of mutually
disjoint bounded intervals. Note that bounded intervals are themselves elementary sets.

A finite set $\balfa=\{\alpha_0,\alpha_1,\dots,\alpha_m\}\subset\ab$ with $m\in\N$
is said to be a {\em division} of the interval $\ab$ if
\begin{align*}
   a=\alpha_0<\alpha_1<\ldots<\alpha_m=b.
\end{align*}
The set of all divisions of $\ab$ is denoted by $\mathcal{D}\ab.$ The symbol $\nu(\balfa)$
will be kept for the number of subintervals $[\alpha_{j{-}1},\alpha_j]\,$
generated by the division $\balfa,$ i.e., $\nu(\balfa)=m$ in the above case.

\smallskip

Let $f {:}\,\ab\to X$ be a function with values in a Banach space $X.$ As in the
case of the real-valued functions, the {\em variation} of $f$ on $\ab$ is defined
by
\[
 \var_a^b f
 =\sup_{\balfa\in\mathcal{D}\ab}
        \sum_{j{=}1}^{\nu(\balfa)}\|f(\alpha_j)-f(\alpha_{j{-}1})\|_X.
\]
If $\var_a^b f<\infty,$ then $f$ has a bounded variation on $\ab$.
$BV(\ab,X)$ is the set of all functions $f {:}\,\ab\to X$ of a bounded variation
on $\ab.$

\smallskip

Let $\mathcal{B}\,{=}\,(X,Y,Z)$ be a bilinear triple. For $f{:}\,\ab\to X$ and
a division $\balfa\,{=}\,\{\alpha_0,\alpha_1,\Big.$ $\Big. \ldots,\alpha_m\}$ of $\ab,$
we define
\begin{equation*}
      \left.\begin{array}{ll}\displaystyle
      \,&(\mathcal{B})V_{a}^{b}(f,\balfa)
      \\[0mm]
      &\displaystyle\quad:=
      \sup\Big\{\Big\|
      \sum_{j{=}1}^{\nu(\balfa)}[f(\alpha_j)\,{-}\,f(\alpha_{j{-}1})]\,y_j
      \Big\|_{Z}\hskip-1mm{:}\,
      y_j\in Y,\,\|y_j\|_{Y}\le 1,\,j=1,2,\dots,\nu(\balfa)
      \Big\}
\end{array}\right\}\hskip-1mm
\end{equation*}
and
\begin{align*}
   (\mathcal{B})\var_{a}^{b}f
   =\sup\Big\{(\mathcal{B})V_{a}^{b}(f,\balfa)\,{:}\,\balfa\in\mathcal{D}\ab\Big\}.
\end{align*}
A function $f{:}\,\ab\to X$ with $(\mathcal{B})\var_{a}^{b}(f)<\infty$
is said to have a {\em bounded $\mathcal{B}-$variation} on $\ab$ or
a {\em bounded semi-variation}. The set of all functions $f{:}\,\ab\to X$ with
bounded $\mathcal{B}-$variation on $\ab$ is denoted by $(\mathcal{B})BV(\ab,X).$

\smallskip

$G(\ab,X)$ denotes the set of all $X$-valued functions which are regulated on $\ab.$
Recall that $f{:}\,\ab\to X$ is {\em regulated} on $\ab$ if for each $t\in[a,b)$
there is a $f(t+)\in X,$ such that
\[
   \lim_{s{\to} t+}\|f(s)-f(t+)\|_X=0,
\]
and for each $t\in (a,b\,],$ there is a $f(t-)\in X,$ such that
\[
   \lim_{s{\to} t-}\|f(s)- f(t-)\|_X=0\,.
\]
For $f\in G(\ab,X)$ and $t\in\ab,$ we denote $\Delta^+f(t)=f(t+)-f(t),$
$\Delta^-f(t)=f(t)-f(t-)$ and $\Delta f(t)=f(t+)-f(t-)$
(where by convention $\Delta^-f(a)=\Delta^+f(b)=0$).

\smallskip

Let $\mathcal{B}=(X,Y,Z)$ be a bilinear triple. A function $f{:}\,\ab\to X$ is called
{\em $\mathcal{B}-$regulated} on $\ab$ (or {\em simply-regulated} on $\ab$) if
the function $fy:t\in\ab\to f(t)y\in Z$ is regulated for all $y\in Y.$ The set of all
simply-regulated functions $f{:}\,\ab\to X$ is denoted by $(\mathcal{B})G(\ab,X).$

Clearly, $G(\ab,X)\subset(\mathcal{B})G(\ab,X)$. Moreover,
\[
   BV(\ab,X)\subset G(\ab,X)\mbox{ \ and \ } BV(\ab,X)\subset(\mathcal{B})BV(\ab,X).
\]

\smallskip

A finite set of points in $\ab$
\begin{align*}
  P=\{
      \alpha_0,\xi_1,\alpha_1,\xi_2,\dots,\alpha_{m-1},\xi_m,\alpha_m
    \}
\end{align*}
where $\{\alpha_0,\alpha_1,\dots,\alpha_m\}\in\mathcal{D}\ab$ and
$\xi_j\in[\alpha_{j{-}1},\alpha_j]\mbox{ \ for\ }j=1,2,\dots,\nu(P)$ is called
a {\em tagged partition} of $\ab.$ The point $\xi_j$ is called the \emph{tag} of the subinterval $[\alpha_{j{-}1},\alpha_j]$ for every $j=1,2,\dots,\nu(P).$
We then shall write
\[
    P=\{([\alpha_{j{-}1},\alpha_j],\xi_j)\}\,\,\,\text{or}\,\,\,P=(\balfa,\ksi)
\]
with $\balfa=\{\alpha_0,\alpha_1,\,\dots\,,\alpha_m\},$
$\ksi=\{\xi_1,\xi_1,\,\dots\,,\xi_m\}$ and $\nu(P)=\nu(\balfa).$

\smallskip

Positive functions $\delta\,{:}\,\ab\to (0,\infty)$ are called {\em gauges}
on $\ab.$ For a given gauge $\delta$ on $\ab,$ a tagged partition
$P=\{([\alpha_{j{-}1},\alpha_j],\xi_j)\}$ of $\ab$ is called $\delta-$fine if
\[
 [\alpha_{j{-}1},\alpha_j]\,{\subset}\,(\xi_j{-}\delta(\_j),\xi_j{+}\delta(\xi_j))
 \quad\mbox{for \ } j\,{=}\,1,2,\dots,\nu(P).
\]

The following lemma shows that the set of $\delta-$fine partitions is nonempty and this result is known as the Cousin lemma (see e.g., \cite{Co}, \cite[Lemma 9.2]{G2}, \cite[Theorem 2.3.1]{L1}, \cite[Theorem 1.1.5]{L2}).\smallskip

\begin{lemma}\label{Cousin}{\emph{[}\bf{Cousin}]}
Given an arbitrary gauge $\delta$ on $\ab,$ there is a $\delta-$fine partition of $\ab.$
\end{lemma}\smallskip

If $\mathcal{B}=(X,Y,Z)$ is a bilinear triple, then for functions $f{:}\,\ab\to X,$
$g{:}\,\ab\,{\to}\,Y$ and a tagged partition $P=\{([\alpha_{j{-}1},\alpha_j],\xi_j)\}$
of $\ab,$ we set
\begin{align*}
    S(\dd f,g,P)=\sum_{j{=}1}^{\nu(P)}[f(\alpha_j)-f(\alpha_{j{-}1})]\,g(\xi_j)
\end{align*}
and
\begin{align*}
    S(f,\dd g,P)=\sum_{j{=}1}^{\nu(P)}f(\xi_j)[g(\alpha_j)-g(\alpha_{j{-}1})].
\end{align*}

Now, we can present the definition of the abstract Kurzweil-Stieltjes integral
as introduced by \v{S}.~Schwabik in \cite[Definition 5]{S}.\smallskip

\begin{definition}\label{KS}{\rm
Let $\mathcal{B}=(X,Y,Z)$ be a bilinear triple and let $f{:}\,\ab\to X$ and
$g{:}\,\ab\,{\to}\,Y$ be given. We say that the Kurzweil-Stieltjes integral
(shortly KS-integral) $\int_a^b\f\,g$ exists if there is $I\in Z$ such that
for every $\eps\,{>}\,0$ there is a gauge $\delta$ on $\ab$ such that
\begin{equation}\label{Def_KS_1}
  \Big\|S(\dd f,g,P)-I\Big\|_Z<\eps
\end{equation}
holds for every $\delta-$fine partition $P$ of $\ab.$ In such a case, we put
\[
   \int_a^b\f\,g=I.
\]
Furthermore, we define
\[
   \int_a^a\f\,g=0  \quad\mbox{and}\quad
   \int_a^b\f\,g=-\int_b^a\f\,g \quad\mbox{if \ } b<a.
\]
\skipaline
Similarly, if $f{:}\,\ab\to X$ and $g{:}\,\ab\to Y,$ then $\int_a^b f\,\g=I\in Z$
if and only if for every $\eps\,{>}\,0$ there is a gauge $\delta$ on $\ab$ such that
\begin{equation}\label{Def_KS_2}
  \Big\|S(f,\dd g,P)-I\Big\|_Z<\eps
\end{equation}
holds for every $\delta-$fine partition $P$ of $\ab.$
}\end{definition}

\skipaline

The Kurzweil-Stieltjes integral is well defined by Definition \ref{KS} owing to the Cousin lemma \ref{Cousin}.
The existence of the Kurzweil-Stieltjes integral is guaranteed, (see e.g., \cite[Proposition 15]{S}).
Evidently, it reduces to the Kurzweil-Henstock integral whenever the integrator $f$ in \eqref{Def_KS_1} (the integrator $g$ in \eqref{Def_KS_2}) is an identity function.

\skipaline

Throughout the paper, we assume that $\mathcal{B}=(X,Y,Z)$ is a bilinear triple.
Furthermore, $\ab$ is a fixed bounded and closed interval in $\R.$
All functions $f$ are supposed to be defined on the entire interval $\ab$ and
extended outside the interval $\ab$ in such a way that $f(t)=f(a)$ and $f(s)=f(b)$
for $t<a$ and $s>b.$\\

\section{Integration over arbitrary bounded sets}\label{IOABS}

\skiphalfaline

In \cite[Section 5]{MHT}, the Kurzweil-Stieltjes integral of operator-valued functions
over elementary subsets of $\ab$ was introduced, and its basic properties were described.
This definition can be easily extended to the arbitrary subsets of $\ab$
and to setting in a general bilinear triple $\mathcal{B}=(X,Y,Z)$.\smallskip

\begin{definition}\label{def.int.bounded}{\rm
Let $f{:}\,\ab\to X,\,\,g{:}\,\ab\to Y$ and let $S$ be an arbitrary subset
of $\ab.$ Then, the Kurzweil-Stieltjes integral (shortly KS-integral or integral) of $g$ with respect
to $f$ over the set $S,$ denoted by $\int_{S}\f\,g,$ is defined by
\begin{align*}
    \int_{S}\f\,g:=\int_a^b\f\,(g\chi_{S})
\end{align*}
whenever the integral on the right-hand side exists.

\skipaline

Similarly, if $f\,{:}\,\ab\to X,\,\,g\,{:}\,\ab\to Y,$ then the integral
$\int_{S}f\,\g$ is defined by
\begin{align*}
    \int_{S}f\,\g:=\int_a^b(f\chi_{S})\,\g
\end{align*}
whenever the integral on the right-hand side exists.
}\end{definition}

\skipaline

\begin{remark}\label{R1}{\rm
By Definitions \ref{KS} and \ref{def.int.bounded}, the existence of the integral
$\int_{S}\f\,g$ means that there exists $I\in Z$ with the following property:
for every $\eps\,{>}\,0$ there exists a gauge $\delta$ on $\ab$ such that
\begin{align*}
\Big\|S(\dd f,g\chi_{S},P)-I\Big\|_Z<\eps
\end{align*}
whenever $P$ is a $\delta-$fine partition of $\ab.$
}\end{remark}

\skipaline

Definition 5.1 from \cite{MHT} is a special case of Definition~\ref{def.int.bounded}.
However, all the results presented in \cite{MHT} for the special case
$\mathcal{B}{=}\,(\mathcal{L}(X,Z),X,Z)$  can be reformulated for the setting of this
paper with a general bilinear triple $\mathcal{B}=(X,Y,Z).$ In particular,
Propositions ~\ref{linearity} and \ref{int.interval} are valid.

\skipaline

\begin{proposition}\label{linearity}
Let $S$ be an arbitrary subset of $\ab.$ Then\emph{,} the following assertions are true\emph{:}

\begin{enumerate}\leftskip=-2pt
\item[(i)]
Let $f{:}\,\ab\to X$ and $g_i{:}\,\ab\to Y,$ $i=1,2,$ such that
the integrals $\int_{S}\f\,g_i$ for $i=1,2$ exist. Then\emph{,} the integral
$\int_{S}\f (c_1 g_1+c_2 g_2)$ also exists\emph{,} and
\[
  \int_{S}\f\,(c_1g_1+c_2g_2)
  =c_1\int_{S}\f\,g_1+c_2\int_{S}\f\,g_2
\]
for all $c_1,c_2\in\R.$\vskip2mm

\item[(ii)]
Let $f_i{:}\,\ab\to X,$ $i=1,2,$ and $g{:}\,\ab\to Y,$ such that
the integrals $\int_{S}[\dd f_i]\,g$ for $i=1,2$ exist. Then\emph{,} the integral
$\int_{S}[\dd\left(c_1 f_1+c_2 f_2\right)]\,g$ also exists\emph{,} and
\[
  \int_{S}[\dd\left(c_1 f_1+c_2 f_2\right)]\,g
  =c_1\int_{S}[\dd f_1]\,g+c_2\int_{S}[\dd f_2]\,g
\]
for all $c_1,\,c_2\in\R.$
\end{enumerate}
\end{proposition}

\skipaline

\begin{remark}\label{R}{\rm
As $g\chi_{\ab}=g$ on $\ab,$ the integral $\int_a^b\f\,g$ exists if and only
if the integral $\int_{\ab}\f\,g$ exists. In such a case, these integrals have
the same value, i.e.,
\begin{equation}\label{int.closed.interval-added}
   \int_{\ab}\f\,g=\int_a^b\f\,g\,.
\end{equation}
Meanwhile,
\[
   g(t)\chi_{(a,b)}(t)-g(t)=
   \begin{cases}
         -g(a) &\mbox{if \ } t=a,
         \\
         \quad 0 &\mbox{if \ } t\in(a,b),
         \\
         -g(b) &\mbox{if \ } t=b.
   \end{cases}
\]
and hence, by \cite[Lemma 12]{S}, we get for an arbitrary $d\in(a,b)$
\begin{align*}
    \int_a^b\f\,(g\chi_{(a,b)}-g)
    &=\int_a^d\f\,(g\chi_{(a,b)}-g)+\int_d^b\f\,(g\chi_{(a,b)}-g)
    \\
    &={-}(\lim_{r{\to}a+}[f(r)\,g(a)]\,{-}\,f(a)\,g(a))
      {-}(f(b)\,g(b)\,{-}\,\lim_{r{\to}b-}[f(r)\,(b)]),
\end{align*}
i.e., the integral $\int_a^b\f\,g$ exists if and only if the integral $\int_{(a,b)}\f\,g$ exists,
and in such a case,
\begin{equation}\label{r3,1}
      \left.\begin{array}{ll}\displaystyle
      \int_{(a,b)}\f\,g&=\displaystyle \int_a^b\f\,g
      \\[4mm]
      &\displaystyle\quad+f(a)\,g(a)-f(b)\,g(b)
      +\lim_{r{\to}b-}[f(r)\,g(b)]-\lim_{r{\to}a+}[f(r)\,g(a)].
\end{array}\right\}\hskip-3mm
\end{equation}
}\end{remark}

\skipaline

The next proposition summarizes the properties of the KS-integral over all possible kinds
of subintervals of $\ab.$ The proofs of its assertions are the easy modifications of those
of \cite[Theorems 5.8, 5.10, and 5.11]{MHT}. The above observations concerning the cases
$c=a$ and/or $d=b$ will be included, considering the convention that functions
$f$ and $g$ are to be considered extended outside of the interval $\ab$
as constant functions on $(-\infty,a]\cup[b,\infty).$

\skipaline

\begin{proposition}\label{int.interval}
Let $f\in(\mathcal{B})G(\ab;X),$ $g{:}\,\ab\to Y,$ and $a\le c<d\le b.$
Then\emph{,} the following assertions are true$:$
\begin{enumerate}\leftskip = -2 pt	
\item[(i)]
The integral $\int_{(c,d)}\f\,g$ exists if and only if the integral
$\int_c^d\f\,g$ exists. In such a case\emph{,}
\begin{equation*}
      \left.\begin{array}{ll}\displaystyle
      \int_{(c,d)}\f\,g
      &=\displaystyle f(c)\,g(c)-\lim_{r{\to}c+}[f(r)\,g(c)]
      +\int_c^d\f\,g
      \\[4mm]
      &\displaystyle\quad-f(d)\,g(d)+\lim_{r{\to}d-}[f(r)\,g(d)].
\end{array}\right\}\hskip-3mm
\end{equation*}

\item[(ii)]
The integral $\int_{[c,d)}\f\,g$ exists if and only if the integral
$\int_c^d\f\,g$ exists. In such a case\emph{,}
\begin{equation*}
\left.\begin{array}{ll}\displaystyle
  \int_{[c,d)}\f\,g
   &=\displaystyle f(c)\,g(c)-\lim_{r{\to}c-}[f(r)\,g(c)]+\int_c^d\f\,g
 \\[4mm]
   &\displaystyle\quad-f(d)\,g(d)-\lim_{r{\to}d-}[f(r)\,g(d)].
\end{array}\right\}\hskip-3mm
\end{equation*}

\item[(iii)]
The integral $\int_{(c,d]}\f\,g$ exists if and only if the integral
$\int_c^d\f\,g$ exists. In such a case\emph{,}
\begin{equation*}
\left.\begin{array}{ll}\displaystyle
  \int_{(c,d]}\f\,g
   &=\displaystyle f(c)\,g(c)-\lim_{r{\to}c+}[f(r)\,g(c)]+\int_c^d\f\,g
 \\[4mm]
   &\displaystyle\quad+\lim_{r{\to}d+}[f(r)\,g(d)]-f(d)\,g(d).
\end{array}\right\}\hskip-3mm
\end{equation*}

\item[(iv)]
The integral $\int_{[c,d]}\f\,g$ exists if and only if the integral
$\int_c^d\f\,g$ exists. In such a case\emph{,}
\begin{equation*}
\left.\begin{array}{ll}\displaystyle
            \int_{[c,d]}\f\,g
            &=\displaystyle f(c)\,g(c)-\lim_{r{\to}c-}[f(r)\,g(c)]+\int_c^d\f\,g
            \\[4mm]
            &\displaystyle\quad+\lim_{r{\to}d+}[f(r)\,g(d)]-f(d)\,g(d).
      \end{array}\right\}\hskip-3mm
\end{equation*}
\end{enumerate}
\end{proposition}

\skipaline\vskip3mm

\begin{remark}\label{rem3.5}
{\rm\quad
If $a\le c<d\le b,$ $f\in(\mathcal{B})G(\ab;X),$ and $g{:}\,\ab\to Y,$
then Proposition \ref{int.interval} implies that if any one of the integrals
\begin{equation}\label{differ-equal-1}
      \int_{(c,d)}\f\,g,\,\, \int_{[c,d)}\f\,g,\,\,
      \int_{(c,d]}\f\,g,\,\, \int_{[c,d]}\f\,g,\,\,
      \int_c^d\f\,g
\end{equation}
exists, then all the others exist as well. Of course, their values can differ, generally. If, in addition, $f$ is continuous on $\ab,$ then all the equalities
\begin{equation}\label{differ-equal-2}
      \int_{(c,d)}\f\,g=\int_{[c,d)}\f\,g
      =\int_{(c,d]}\f\,g=\int_{[c,d]}\f\,g
      =\int_c^d\f\,g
\end{equation}
are true.
}
\end{remark}

\skipaline

The existence of the integral
$\int_a^b\f\,g$ need not (even in the case of the identity
integrator $f(x)=x$) always imply the existence of the
integral $\int_{T}\f\,g$ for every subset $T$ of $\ab.$ The
examples for that can be constructed on the basis of the
Lebesgue integral (the special case of the Kurzweil-Henstock
integral) assuming e.g., that $T$ is not measurable. As shown by
the next assertion, this cannot happen when we restrict
ourselves to the elementary subsets of $\ab.$

\skipaline

\begin{theorem}\label{int.elem.set}
The following assertions are true for all $f{\in}\,(\mathcal{B})G(\ab;X)$
and $g{:}\,\ab\,{\to}\,Y.$
\begin{itemize}\leftskip=-2pt	
      \item[(i)]
           Let $E$ be an elementary subset of $\ab$ such that the integral $\int_{E}\f\,g$
           exists. Then\emph{,} the integral $\int_{T}\f\,g$ exists for every elementary subset $T$
           of $E.$\vskip2mm

      \item[(ii)]
           Let $E=\bigcup_{k=1}^{p}J_k,$ where $\{J_k\,{:}\,k=1,2,\dots,p\}$
           are mutually disjoint subintervals of~$\ab,$ and let the integral $\int_{E}\f\,g$
           exist. Then\emph{,} all the integrals
                  \[
                  \int_{J_k}\hskip-1mm \f g,\quad k=1,2,\dots,p,
                  \]
           exist as well and
                  \begin{equation}\label{E}
                  \int_{E}\f\,g=\sum_{k=1}^{p}\int_{J_k}\f\,g.
                  \end{equation}
\end{itemize}
\end{theorem}
\proof \quad (i) See \cite[Corollary 5.15]{MHT}.

\noindent (ii)
By \cite[Theorem 5.13]{MHT}, this assertion is true if the set
$\{J_k\,{:}\,k=1,2,\dots,p\}$ is a~minimal
decomposition of $E,$ i.e., (cf. \cite[Definition 4.9]{MHT}) the union
$J_k\cup J_\ell$ is not an~interval whenever $k\,{\ne}\,\ell.$ Of course,
we may assume that the intervals $\{J_k\}$ are ordered in such a way that
$x\le y$ holds whenever $x\in J_k,$ $y\in J_\ell$ and $k<\ell.$ Then, if
$\{J_k\,{:}\,k=1,2,\dots,p\}$ is not a~minimal
decomposition, there must exist $k\,{\in}\,\{1,2,\dots,p\,{-}\,1\}$
such that $J\,{=}\,J_k\,{\cup}\,J_{k+1}$ is an~interval. Then, as
$J_k\,{\cap}\,J_{k{+}1}\,{=}\,\emptyset,$ we get
   \[
    \int_{J_k}\f\,g+\int_{J_{k+1}}\f\,g
   =\int_a^b\f\,g(\chi_{J_k}+\chi_{J_{k+1}})
   =\int_a^b\f\,(g\chi_J)
   =\int_J\f\,g.
   \]
Hence, when we replace in the sum on the right-hand side
of \eqref{E} all such couples by their unions, this sum does not change,
wherefrom assertion (ii) follows.
\endproof

\skipaline

\begin{remark}\label{min.decomp}
{\rm
Let \ $E=\bigcup_{k=1}^{p}J_k$ be an elementary set of $\ab.$
\begin{enumerate}\leftskip=-5mm	
\item[(i)]   Let \ $\{J^{*}_k{:}\,k=1,2,\dots,p^*\}$ be the minimal decomposition of $E.$
             Then, Theorem~\ref{int.elem.set}~(ii) implies that $p^*\,{\le}\,p$ and
             \[
             \int_{E}\f\,g=\sum_{k=1}^{p}\int_{J_k}\f\,g
             =\sum_{k=1}^{p^*}\int_{J^*_k}\f\,g\,.
             \]
             \vskip2mm
\item[(ii)]  Let $c_k$ and $d_k$ be an infimum and a supremum of $J_k,$ respectively, for every $k=1,2,\dots,p.$ Then, from Remark \ref{rem3.5} with
             Theorem~\ref{int.elem.set}~(ii), the equality
             \begin{equation}\label{rem.inf.sup}
             \int_{E}\f\,g=\sum_{k=1}^{p}\int_{J_k}\f\,g=\sum_{k=1}^{p}\int_{c_k}^{d_k}\f\,g
             \end{equation}
             holds for all continuous integrators $f,$ especially for the Kurzweil-Henstock integral. However, \eqref{rem.inf.sup} is not valid any longer for the KS-integral.
\end{enumerate}
}
\end{remark}

\skipaline

\begin{remark}\label{R3.8}
{\rm
If a function $g{:}\,\ab\to Y$ and a subset $S$ of $\ab$ are such that $g=0$
on $S,$ then $g\chi_S=0$ on $\ab$ and hence
$\int_{S}\f\,g=\int_a^b\f\,(g\chi_S)=0$ and $\int_{T}\f\,g=0$ as well for every subset $T$ of $S$ and every
$f{:}\,\ab\to X.$ In particular, if the integral $\int_{S}\f\,g$
exists and $h{:}\,\ab\to Y$ coincides with $g$ on $S,$ then
$\int_{S}\f\,h=\int_{S}\f\,g.$
}
\end{remark}

\skipaline

The next assertion discloses the additivity properties of the KS-integral over
arbitrary subsets of $\ab.$

\skipaline

\begin{proposition}\label{P1}
Let $f{:}\,\ab\to X,$ $g{:}\,\ab\to Y$ and subsets
$S_1,S_2$ in $\ab$ be given.
Then\emph{,} whenever three of the integrals
\[
  \int_{S_1}\f\,g,\quad\int_{S_2}\f\,g,\quad
  \int_{S_1{\cup}S_2}\f\,g, \quad \int_{S_1{\cap}S_2}\f\,g
\]
exist\emph{,} then there also exists the remaining one and the equality
\[
  \int_{S_1}\f\,g+\int_{S_2}\f\,g
  =\int_{S_1{\cup}S_2}\f\,g+\int_{S_1{\cap}S_2}\f\,g
\]
holds.
\end{proposition}
\proof
It follows directly from the identity
$\,\chi_{S_1}+\chi_{S_2}=\chi_{S_1{\cup}S_2}+\chi_{S_1{\cap}S_2}.$
\endproof

\skiphalfaline

\begin{corollary}\label{C1}
Let $f{:}\,\ab\to X,$ $g{:}\,\ab\to Y,$ and subsets
$S_1,S_2$ in $\ab,$ such that $S_1\,{\cap}\,S_2=\emptyset$ and integrals $\int_{S_1}\f\,g$ and $\int_{S_2}\f\,g$ exist.
Then\emph{,} the integral $\int_{S_1{\cup}S_2}\f\,g$ also exists
and the equality
\[
  \int_{S_1{\cup}S_2}\f\,g=\int_{S_1}\f\,g+\int_{S_2}\f\,g
\]
holds.
\end{corollary}

\skiphalfaline

\begin{proposition}\label{P2berubah}
Let $f{:}\,\ab\to X$ and $g{:}\,\ab\to Y$ be given and let $S_1,\,S_2,\dots,S_p$
be subsets of $\ab$ and $p\ge 2.$ Denote
\[
   S=\bigcup_{j{=}1}^p S_j \quad\mbox{and}\quad
   T_i=\Big(\bigcup_{j{=}1}^{i{-}1}S_j\Big)\cap S_i \mbox{ \ for\ } i=2,3,\dots,p
\]
and assume that all the integrals
\[
  \int_{S_i}\f\,g,\quad \int_{T_i}\f\,g,\quad i=1,2,\dots,p,
\]
exist. Then\emph{,} the integral
\[
  \int_{S}\f\,g
\]
exist as well and
      \begin{equation}\label{additivity}
      \int_{S}\f\,g
      =\sum_{i=1}^{p}\int_{S_i}\f\,g-\sum_{i{=}2}^p\int_{T_i}\f\,g.
      \end{equation}
\end{proposition}
\proof \ (i) \ First, if $p=2,$ then $T_2=S_1\cap S_2,$ and the assertion
of the theorem follows Proposition~\ref{P1}.

\skiphalfaline

\noindent(ii) \  For $p=3,$ we have $S=S_1\cup S_2\cup S_3,$ $T_2=S_1\cap S_2,$
$R_3=(S_1\cup S_2)\cap S_3.$
Denote $M=S_1\cup S_2.$ Then\emph{,} $S=M\cup S_3,$ $T_3=M\cap S_3,$ and by Proposition~\ref{P1},
we get
\begin{align*}
      &\int_{M}\f\,g=\int_{S_1}\f\,g+\int_{S_2}\f\,g-\int_{T_2}\f\,g
      \\[2mm]\noalign{\noindent\mbox{and}}
      &\int_{S}\f\,g
      =\int_{M}\f\,g+\int_{S_3}\f\,g-\int_{M\cap S_3}\f\,g
      \\
      &\hskip16mm=\sum_{i{=}1}^3\int_{S_i}\f\,g-\int_{T_2}\f\,g-\int_{T_3}\f\,g
      \\
      &\hskip16mm=\sum_{i{=}1}^{3}\int_{S_i}\f\,g-\sum_{i{=}2}^3\int_{T_i}\f\,g.
\end{align*}

\skiphalfaline

\noindent(iii) \  Let $N>3$ and let \eqref{additivity} hold for $p=N\,{-}\,1.$
Denote
\[
   M=\bigcup_{j{=}1}^{N{-}1}S_j,\quad S=\bigcup_{j{=}1}^N S_j,\quad\mbox{and}\quad
   T_i=\Big(\bigcup_{j{=}1}^{i{-}1}S_j\Big)\cap S_i \mbox{ \ for\ } i=2,3,\dots,N-1.
\]
Then, $S=M\cup S_N,$ $T_N=M\cap S_N,$ and by \eqref{additivity} (with $n\,{=}\,N\,{-}\,1$)
and Proposition~\ref{P1}, we have
\begin{align*}
      \int_{M}\f\,g
      &=\sum_{i{=}1}^{N{-}1}\int_{S_i}\f\,g-\sum_{i{=}2}^{N{-}1}\int_{T_i}\f\,g
      \\\noalign{\noindent\mbox{and}}
      \int_{S}\f\,g
      &=\int_{M}\f\,g+\int_{S_N}\f\,g-\int_{M\cap S_N}\f\,g
      \\
      &=\sum_{i{=}1}^N\int_{S_i}\f\,g-\sum_{i{=}2}^{N{-}1}\int_{T_i}\f\,g-\int_{T_N}\f\,g
      \\
      &=\sum_{i{=}1}^N\int_{S_i}\f\,g-\sum_{i{=}2}^N\int_{T_i}\f\,g.
\end{align*}
With the induction principle, this completes the proof of the theorem.
\endproof

\skipaline\vskip-0.5mm

\section{Equi-integrability and Equiregulatedness}\label{equi}

\skiphalfaline

Convergence theorems belong to the most important topics discussed in the
frames of integration theory. For the abstract KS-integral, the uniform convergence
theorem given by \v{S}.~Schwabik in \cite[Theorem 11]{S} is the simplest one. It states that
if the sequence $\{g_n\}$ tends uniformly to $g$ on $\ab$ and if all the integrals
$\int_a^b\f\,g_n,$ $n\in\N,$ exist, then the integral $\int_a^b\f\,g$ exists as well
and
\[
    \int_a^b\f\,g=\lim_{n{\to}\infty}\int_a^b\f\,g_n\,.
\]\smallskip
When the uniform convergence of  $\{g_n\}$ to $g$ is replaced by a just pointwise
convergence on $\ab,$ the situation is more difficult. One possible way is indicated
by the bounded convergence theorem (for the abstract KS-integral, see
\cite[Theorem 6.3]{MHT}), which requires the uniform boundedness of the sequence
$\{g_n\}$ on $\ab.$ The next theorem deals with the case that the sequence $\{g_n\}$
even need not be bounded, and its integrator involved is not just a function $f,$
but a sequence $\{f_n\}.$

\skipaline

\begin{theorem}[{Equi-integrability convergence theorem}]\label{H3-convergence}
Let $f_n{:}\,\ab$ ${\to}X$ and $g_n{:}\,\ab{\to}Y,$ for $n\in\N,$ such
that the integral $\int_a^b[\dd f_n]\,g_n$ exists for each $n\in\N.$
Furthermore\emph{,} let the functions $f{:}\,\ab{\to}X$ and $g{:}\,\ab{\to}Y$ such that
the sequences $\{f_n\}$ and $\{g_n\}$ converge pointwise on $\ab$ to $f$ and $g,$
respectively. Finally\emph{,} suppose that
\begin{equation}\label{H3-1}
      \left.\begin{array}{l}\displaystyle
      \mbox{for every $\eta>0$ there is a gauge $\delta$ on $\ab$ such that\ }
      \\[2mm]\displaystyle\quad
      \Big\|S(\dd f_n,g_n,P)-\int_a^b[\dd f_n]\,g_n\Big\|_{Z}<\eta
      \\[4mm]\displaystyle
      \mbox{for every $\delta-$fine partition P of\ } \ab
      \mbox{\ and every\ } n\in\N\,.
      \end{array}\right\}
\end{equation}
Then\emph{,} the integrals $\int_a^b\f\,g$ and $\lim_{n{\to}\infty}\int_a^b[\dd f_n]\,g_n$
exist and
\begin{equation}\label{H3-2}
      \int_a^b\f\,g=\lim_{n{\to}\infty}\int_a^b[\dd f_n]\,g_n.
\end{equation}
\end{theorem}
\begin{proof} \ {\em Step}~1. \ Let $\eps\,{>}\,0$ be given and let $\delta$ be
the gauge corresponding to $\eta=\frac{\eps}{4}$ by \eqref{H3-1}. Then,
\begin{equation}\label{step.first}
      \left.\begin{array}{l}\displaystyle
      \quad\Big\|S(\dd f_n,g_n,P)-\int_a^b[\dd f_n]\,g_n\Big\|_{Z}<\frac{\eps}{4}
      \\[4.5mm]\displaystyle
      \mbox{for every\ } n\in\N \mbox{\ and every\ } \delta-\mbox{fine partition\ }
      P \mbox{\ of\ }\ab.
      \end{array}\right\}\hskip7mm
\end{equation}
By the Cousin lemma \ref{Cousin}, we may fix an arbitrary $\delta-$fine partition $\wt{P}$
of $\ab.$ Due to the pointwise convergence on $\ab$ of $\{f_n\}$ to $f$
and of $\{g_n\}$ to $g,$ we have
\[
      \lim_{n\to\infty}S(\dd f_n,g_n,\wt{P})=S(\dd f,g,\wt{P}).
\]
Hence, we can choose $n_0\in\N,$ such that the inequality
\begin{equation}\label{step}
      \Big\|S(\dd f_n,g_n,\wt{P})-S(\dd f,g,\wt{P})\Big\|_Z<\frac{\eps}{4}
\end{equation}
holds for all $n\ge n_0.$ Let arbitrary $n_1,n_2\ge n_0$ be given. Then,
using \eqref{step.first} and \eqref{step}, we deduce that
\begin{align*}
      \left.\begin{array}{ll}\displaystyle
      \,&\Big\|\int_a^b[\dd f_{n_1}]\,g_{n_1}-\int_a^b[\dd f_{n_2}]\,g_{n_2}\Big\|_{Z}
      \\[5mm]
      &\hskip2mm\le\Big\|\int_a^b[\dd f_{n_1}]\,g_{n_1}-S(\dd f_{n_1},g_{n_1},\wt{P})\Big\|_{Z}
      +\Big\|S(\dd f_{n_1},g_{n_1},\wt{P})-S(\dd f,g,\wt{P})\Big\|_{Z}
      \\[5mm]
      &\hskip4mm+\Big\|S(\dd f,g,\wt{P})-S(\dd f_{n_2},g_{n_2},\wt{P})\Big\|_{Z}
      +\Big\|S(\dd f_{n_2},g_{n_2},\wt{P})-\int_a^b[\dd f_{n_2}]\,g_{n_2}\Big\|_{Z}<\eps,
      \end{array}\right\}
\end{align*}
which shows that $\Big\{\int_a^b[\dd f_n]\,g_n\Big\}$ is a Cauchy sequence in $Z.$
Let
\begin{equation*}\label{lim}
      I=\lim_{n\to\infty}\int_a^b[\dd f_n]\,g_n\,.
\end{equation*}

\noindent{\em Step}~2. \ We shall prove that $\int_a^b\f g=I.$ Let $\eps\,{>}\,0$
be given and let $\delta$ be a corresponding gauge given by \eqref{H3-1}. Furthermore,
let $n_{\eps}\in\N$ such that
\[
    \Big\|\int_a^b[\dd f_n]\,g_n-I\Big\|_Z<\eps\quad\mbox{for all\ } n\ge n_{\eps}\,.
\]
Now, let $\wt{P}$ be a fixed $\delta$-fine partition of $\ab.$ Owing to
the pointwise convergence on $\ab$ of $\{f_n\}$ to $f$ and of $\{g_n\}$
to $g,$ we can choose $k_0\ge n_{\eps},$ such that
\[
    \Big\|S(\dd f_{k_0},g_{k_0},\wt{P})-S(\dd f,g,\wt{P})\Big\|_{Z}<\eps.
\]
Hence,
\begin{align*}
      \left.\begin{array}{ll}\displaystyle
            \,&\Big\|S(\dd f,g,\wt{P})-I\Big\|_Z
            \le\Big\|S(\dd f,g,\wt{P})-S(\dd f_{k_0},g_{{\ell}_0},\wt{P})\Big\|_Z
            \\[4mm]
            &\,\,\quad+\Big\|S(\dd f_{k_0},g_{{\ell}_0},\wt{P})-\int_a^b[\dd f_{k_0}]\,g_{k_0}\Big\|_Z
            +\Big\|\int_a^b[\dd f_{k_0}]\,g_{{\ell}_0}-I\Big\|_Z<3\,\eps.
      \end{array}\right\}
\end{align*}
It follows that $\int_a^b\f\,g=I,$ and this completes the proof.
\end{proof}\vskip-5mm

\begin{definition}\label{H4-Def1}{\rm
Let $f,f_n\,{:}\,\ab\to X$ and $g,g_n\,{:}\,\ab\to Y$ for $n\in\N.$ The sequence
$\{g_n\}$ is said to be {\em equi-integrable with respect to $\{f_n\}$ on $\ab$}
\ if the integrals $\int_a^b[\dd f_n]\,g_n$ exist for all $n\in\N$ and the condition
\eqref{H3-1} in Theorem~$\ref{H3-convergence}$ is satisfied.\\[3mm]
Similarly, for any $E\subset\ab,$ the sequence $\{g_n\}$ is
\noindent{\em equi-integrable with respect to $\{f_n\}$ on $E$} if $\{g_n\,\chi_E\}$
is equi-integrable with respect to $\{f_n\}$ on $\ab.$\\[3mm]
}\end{definition}\vskip-5mm

In view of Definition \ref{H4-Def1}, Theorem~\ref{H3-convergence} may be reformulated
as follows:

\skiphalfaline

\noindent{\em Let $f,f_n\,{:}\,\ab\to X$ and $g,g_n\,{:}\,\ab\to Y,$ $n\in\N,$ such that
$\lim_{n{\to}\infty}f_n(t)=f(t)$ and $\lim_{n{\to}\infty}g_n(t)=g(t)$ on $\ab$ and
suppose that the sequence $\{g_n\}$ is equi-integrable with respect to $\{f_n\}$ on
$\ab.$ Then\emph{,} the integral $\int_a^b\f\,g$ and $\lim_{n\to\infty}\int_a^b[\dd f_n]\,g_n$
exist and} \eqref{H3-2} \emph{holds}.

\skiphalfaline

\begin{remark}\label{uniformly-integrable}\quad
{\rm\begin{enumerate}
\item[(i)]    Equi-integrability is often met in the literature dealing with the theory of Kurzweil-Henstock integrals, (see e.g.,
              R. G. Bartle \cite[Chapter 8]{B2}, R. A. Gordon \cite[Chapter 13]{G2} and \cite{G3}, J. Kurzweil \cite[Chapter 5]{K2}, J. Kurzweil and
              J. Jarn\'ik \cite{K5J}, \v{S}. Schwabik \cite[Chapter 1]{S3}, \v{S}. Schwabik and I. Vrko\v{c} \cite{S4V}, \v{S}. Schwabik and G. Ye
              \cite[Chapter 3]{S5Y}).
              Nonetheless, little is known about the conditions that ensure the equi-integrability for Stieltjes-type integrals for real-valued functions, (see \cite{BoPi}, \cite[Chapter 6]{M3}, \cite{M4}).
\item[(ii)]   Referring to, e.g., R. A. Gordon \cite[Definition 13.15]{G2} or \v{S}. Schwabik and I. Vrko\v{c} \cite[Remark 6]{S4V}, a sequence $\{g_n\}$
              equi-integrable with respect to $\{f_n\}$ on $\ab$ can be also called {\em uniformly integrable} with respect to $\{f_n\}$ on $\ab.$
\item[(iii)]  If $f_n=f$ for all $n\in\N,$ then $\{g_n\}$ is equi-integrable with respect to $f$ on $\ab.$ \\[-3mm]
\end{enumerate}
}
\end{remark}

\skiphalfaline

In general, it is rather difficult to verify that the condition \eqref{H3-1}
is satisfied. The following statement at least enables us to decide whether
a given sequence $\{g_n\}$ is equi-integrable with respect to $\{f_n\}$ on $\ab$
without evaluating the values of all the integrals $\int_a^b[\dd f_n]\,g_n,$ $n\in\N.$

\skipaline

\begin{theorem}[Cauchy equi-integrability criterion]\label{H4-Cauchy}
Let $f_n\,{:}\,\ab\to X$ and $g_n\,{:}\,\ab$ $\to Y,$ for $n\in\N.$ Then\emph{,} the sequence
$\{g_n\}$ is equi-integrable with respect to $\{f_n\}$ on $\ab$ if and only if
\begin{equation}\label{H4-Cauchy-1}
      \left.\begin{array}{l}\displaystyle
           \mbox{for every\ }\eps\,{>}\,0 \mbox{\ there is a gauge\ } \delta \mbox{\ such that \ }
           \\[2mm]\displaystyle\quad
           \Big\|S(\dd f_n,g_n,P)-S(\dd f_n,g_n,Q)\Big\|_Z<\eps
           \\[4mm]\displaystyle
           \mbox{holds for all \ $n\in\N$ and all \ $\delta$-fine partitions\ }
           P \mbox{\ and\ } Q \mbox{\ of\ } \ab.
      \end{array}\right\}
\end{equation}
\end{theorem}
\begin{proof} \ (i) \ Assume that the sequence $\{g_n\}$ is equi-integrable
with respect to $\{f_n\}$ on $\ab.$ Let $\eps\,{>}\,0$ be given and let $\delta$
be an arbitrary gauge corresponding by \eqref{H3-1} to $\eta=\eps/2.$  Then, for
any couple $P, Q$ of the $\delta$-fine partitions of $\ab$ and any $n\in\N,$
we obtain
\begin{align*}
      \left.\begin{array}{ll}\displaystyle
            \,&\Big\|S(\dd f_n,g_n,P) - S(\dd f_n,g_n,Q)\Big\|_{Z}
            \\[4mm]
            &\,\,\quad\le\Big\|S(\dd f_n,g_n,P) - \int_a^b[\dd f_n]\,g_n\Big\|_{Z}
            +\Big\|S(\dd f_n,g_n,Q) - \int_a^b[\dd f_n]\,g_n\Big\|_{Z}<\eps\,.
      \end{array}\right\}
\end{align*}
\smallskip

\noindent(ii)\ Assume that condition \eqref{H4-Cauchy-1} is satisfied. Then, by
the Cauchy-Bolzano criterion \cite[Proposition 7]{S} for the existence of the KS-integral,
the integral $\int_a^b[\dd f_n]\,g_n$ exists for every $n\in\N.$
For a given $n\in\N,$ gauge $\delta$ and $\eps>0,$ denote
\[
   \mathcal{I}_{n}(\eps,\delta)
   =\{S(\dd f_n,g_n,P):\,P\mbox{\ is a\ } \delta{-}\mbox{fine tagged partition of}\ \ab\}.
\]
Due to \eqref{H4-Cauchy-1}, we have
\begin{equation}\label{Diam}\hskip-1mm
      \left.\hskip-1mm\begin{array}{l}
            \mbox{diam}(\mathcal{I}_{n}(\eps,\delta))
            \,{=}\sup\big\{\|S(\dd f_n,g_n,P){-}S(\dd f_n,g_n,Q)\|_Z{:}
            \\[2mm]
            \hskip38mm\,P,Q \mbox{\ are\ }\delta{-}\mbox{fine partitions of\ }\ab\big\}{<}\,\eps.
            \hskip-2mm
      \end{array}\right\}
\end{equation}
By the Cousin lemma \ref{Cousin}, any $\mathcal{I}_{n}(\eps,\delta)$ is nonempty and, furthermore,
\[
    0<\eps_1<\eps_2\implies
    \mathcal{I}_{n}(\eps_1,\delta)\subset\mathcal{I}_{n}(\eps_2,\delta)
    \quad\mbox{for every\ } n\in\N \mbox{\ and gauge\ } \delta.
\]
Thus, using the Cantor intersection theorem for complete metric spaces (see e.g.,
\cite[Theorem 5.1.17]{WFT}), we conclude that, for every $n\in\N,$ the intersection
$\displaystyle\bigcap_{\eps>0}\overline{\mathcal{I}_{n}(\eps,\delta)}$
is a one-point
set $\{I_{n}\}$ with
\[
    I_{n}=\int_a^b[\dd f_n]\,g_n\in Z.
\]

As a consequence, if an arbitrary $\eta>0$ and a gauge $\delta_\eps$ are given, such
that (\ref{H4-Cauchy-1}) is true with $\eps=\eta/2,$ then
$I_n\in\overline{\mathcal{I}_{n}(\eps,\delta_\eps)}$
for every $n\in\N.$ In particular,
considering \eqref{Diam}, we have
\begin{align*}
      &\Big\|S(\dd f_n,g_n,P)-\int_a^b[\dd f_n]\,g_n\Big\|_Z
      =\|S(\dd f_n,g_n,P)-I_{n}\|_Z\le\eps<\eta
\end{align*}
for every $\delta_{\eps}-$fine partition $P$ of $\ab$ and every $n\in\N.$ In other words,
the sequence $\{g_n\}$ is equi-integrable with respect to $\{f_n\}$ on $\ab.$
This completes the proof.
\end{proof}

\skipaline

If $E=[c,d]$ is a closed subinterval of $\ab,$ it seems to be natural to also define
the equi-integrability on $E$ in the following alternative way$:$

$\{g_n\}$ is equi-integrable with respect to $\{f_n\}$ on $[c,d]$ if the conditions
\begin{equation}\label{H3-1*}
      \left.\begin{array}{l}\displaystyle
            \mbox{\hskip12mm\quad the integrals \ } \int_c^d[\dd f_n]\,g_n
            \mbox{ \ exist for all \ } n\in\N
            \\
            \mbox{and}
            \\[1.5mm]
            \mbox{\hskip12mm\quad for every $\eta>0$ there is a gauge $\delta$ on $\ab$ such that\ }
            \\[2mm]\displaystyle\qquad
            \hskip25mm\Big\|S(\dd f_n,g_n,P)-\int_c^d[\dd f_n]\,g_n\Big\|_{Z}<\eta
            \\[4mm]\displaystyle
            \mbox{\hskip12mm\quad holds for every $\delta-$fine partition $P$ \ of\ } [c,d]
            \mbox{\ and every\ } n\in\N\,\hskip2mm
      \end{array}\right\}
\end{equation}
are satisfied.

\skiphalfaline

Thus, any comparison of the two possible definitions is urgently needed. To this
aim, we shall first recall the notion of equiregulatedness based on
Fra\v{n}\-kov\'{a}~\cite{Fra}.

\skiphalfaline

\begin{definition}\label{equireg}
{\em
A subset $M$ of~$G(\ab;X)$ is called {\em equiregulated} if the following
conditions hold.
\begin{enumerate}
\item[(i)]   For each $\eps>0$ and $\tau\in(a,b]$ there is a~$\delta_1(\tau)\in(0,\tau-a)$
             such that
             \[
             \|f(\tau-)-f(t)\|_X<\eps
             \quad\mbox{for all\ } t\in(\tau-\delta_1(\tau),\tau) \mbox{ \ and \ } f\in M.
             \]
\item[(ii)]  For each $\eps>0$ and $\tau\in[a,b)$ there is a~$\delta_2(\tau)\in(0,b-\tau)$
             such that
             \[
             \|f(\tau+)-f(t)\|_X<\eps
             \quad\mbox{for all\ } t\in(\tau,\tau+\delta_2(\tau)) \mbox{ \ and \ } f\in M.
             \]
\end{enumerate}}
\end{definition}

\skipaline

Definition \ref{equireg} allows us to develop the concept of equi-integrability of sequence
$\{g_n\}$ of functions mapping $\ab$ into $Y$ with respect to $\{f_n\}$ of functions mapping $\ab$ into $X$ over an arbitrary elementary set $E$ of $\ab$.

\begin{lemma}\label{equi.properties-2A}
Let $J=[c,d]\subset\ab,$ let the sequence $\{f_n\}\subset G(\ab;X)$ be equiregulated
and let the sequence $\{g_n\}$ of functions mapping $\ab$ into $Y$ be pointwise bounded
on $[c,d],$
i.e.\emph{,} the sequence $\{g_n(t)\}$ is bounded in $Y$ for each $t\in[c,d].$  Then\emph{,} $\{g_n\}$
is equi-integrable with respect to $\{f_n\}$ on $J$ if and only if condition
\eqref{H3-1*} is satisfied.
\end{lemma}
\proof\ (i) \ Let the assumptions of the lemma be satisfied and let \eqref{H3-1*}
be true. First, assume that $a<c<d<b.$

Let an arbitrary $\eps>0$ be given and let $\wt{\delta}$ be a gauge on $[c,d]$
such that \eqref{H3-1*} holds for $\eta=\eps/3,$ i.e.,
\begin{equation}\label{equivalence-1}
      \left.\begin{array}{l}\displaystyle
            \quad\Big\|S(\dd f_n,g_n,\wt{P})-\int_c^d[\dd f_n]\,g_n\Big\|_{Z}<\frac{\eps}3
            \\[5mm]
            \mbox{for all\ } n\in\N \mbox{\ and all\ } \wt{\delta}-\mbox{fine partitions\ }
            \wt{P} \mbox{\ of\ }[c,d]\,.
      \end{array}\right\}
\end{equation}
Put
\[
 \delta(t)=\begin{cases}
                 \quad \min\{\frac 14(c-t),1\} &\mbox{if \ } t\in[a,c)\,,
                 \\[2mm]
                 \quad \min\{\Delta,\wt{\delta}(c)\} &\mbox{if \ } t=c\,,
                 \\[2mm]
                 \quad \min\{\frac 14(t-c),\frac 14(d-t),\wt{\delta}(t)\}
                 &\mbox{if \ } t\in(c,d)\,,
                 \\[2mm]
                 \quad \min\{\Delta,\wt{\delta}(d)\} &\mbox{if \ } t=d\,,
                 \\[2mm]
                 \quad \min\{\frac 14(t-d),1\} &\mbox{if \ } t\in(d,b]\,,
          \end{cases}
\]\vskip0.6mm
\noindent where $\Delta\in(0,\min\{c-a,b-d\}),$ such that
\begin{equation}\label{equivalence-2}
      \left.\begin{array}{l}\hskip10mm\displaystyle
            \|f_n(t)-f_n(c-)\|_{X}\,{<}\,\frac{\eps}{3\,(1{+}\sup\{\|g_n(c)\|_Y\,{:}\,n\in\N\})}
            \\[5mm]
            \hskip48mm\mbox{\ for\ } n\in\N\mbox{\ and\ } t\in(c\,{-}\,\Delta,c),
            \\[1mm]\mbox{and}\quad
            \\[2mm]\hskip10mm\displaystyle
            \|f_n(d+)-f_n(t)\|_{X}\,{<}\,\frac{\eps}{3\,(1{+}\sup\{\|g_n(d)\|_Y\,{:}\,n\in\N\})}
            \\[5mm]
            \hskip48mm\mbox{\ for\ } n\in\N \mbox{\ and\ } t\in(d,d\,{+}\,\Delta)\,.
      \end{array}\right\}
\end{equation}\vskip0.6mm
\noindent Such a $\Delta$ may be chosen because we assume the equiregulatedness of $\{f_n\}$
and the boundedness of $\{g_n(c)\}$ and $\{g_n(d)\}$ in $Y.$

\skiphalfaline

Now, let $P=\big\{([\alpha_{j{-}1},\alpha_j],\xi_j)\big\}$ be a $\delta$-fine
partition of $\ab$ such that $\nu(P)>4.$ Then, cf. \cite[Lemma 6.2.11]{M3},
there are indices $k,\ell\in\N$ such that $2\le k<\ell\le\nu(P)-2,$
$\xi_{k}=c>\alpha_{k{-}1}>c-\Delta>a$ and $\xi_{\ell}=d<\alpha_{\ell{+}1}<d+\Delta<b.$
In addition, we may assume that $\xi_{k{-}1}=c=\alpha_k$ and $d=\alpha_\ell=\xi_{\ell{+}1}.$
To summarize, we have
\begin{equation*}\label{P}
      \left.\begin{array}{l}
            a<c-\Delta<\alpha_{k{-}1}<\xi_k=c=\alpha_k=\xi_{k{+}1}<\alpha_{k{+}1}
            \\[2mm]\quad
            <\xi_{\ell}=d=\alpha_{\ell}=\xi_{\ell{+}1}<\alpha_{\ell{+}1}<d+\Delta<b.
      \end{array}\right\}
\end{equation*}
Hence, we can successively deduce
\begin{equation*}
      \left.\begin{array}{ll}\displaystyle
            &S(\dd f_n,g_n\chi_{[c,d]},P)
            =\sum_{j=k}^{\ell+1}[f_n(\alpha_j)-f_n(\alpha_{j-1})]\,g_n(\xi_j)
            \\[4mm]&\hskip13mm
            =\sum_{j=k+1}^{\ell}[f_n(\alpha_j)-f_n(\alpha_{j-1})]\,g_n(\xi_j)
            \\[4mm]&\hskip20mm
            +[f_n(c)-f_n(\alpha_{k-1})]\,g_n(c)+[f_n(\alpha_{\ell+1})-f_n(d)]\,g_n(d)
            \\[4mm]\displaybreak[0]&\hskip13mm
            =\sum_{j=k+1}^{\ell}[f_n(\alpha_j)-f_n(\alpha_{j-1})]\,g_n(\xi_j)
            +\Delta^-f_n(c)\,g_n(c)+\Delta^+f_n(d)\,g_n(d)
            \\[4mm]&\hskip20mm
            +[f_n(c-)-f_n(\alpha_{k{-}1})]\,g_n(c)+[f_n(\alpha_{\ell+1})-f_n(d+)]\,g_n(d),
      \end{array}\right\}
\end{equation*}
i.e.,\vskip-4mm
\begin{equation}\label{S1}
      \left.\begin{array}{ll}\displaystyle
      &S(\dd f_n,g_n\chi_{[c,d]},P)
      =\displaystyle\sum_{j=k+1}^{\ell}[f_n(\alpha_j){-}f_n(\alpha_{j-1})]\,g_n(\xi_j)
      \\[7mm]\displaystyle
      &\hskip13mm+\Delta^-f_n(c)\,g_n(c)\,{+}\,\Delta^+f_n(d)\,g_n(d)
      \\[4mm]\displaystyle
      &\hskip13mm+[f_n(c-){-}f_n(\alpha_{k{-}1})]\,g_n(c)+[f_n(\alpha_{\ell+1}){-}f_n(d+)]\,g_n(d).
      \end{array}\right\}
\end{equation}\\
As for every $n\in\N$ the integral $\int_c^d[\dd f_n]\,g_n$ exists, by
Proposition~\ref{int.interval}~(iii) (see also Remark~\ref{rem3.5}),
all the integrals $\int_{[c,d]}[\dd f_n]\,g_n,$ $n\in\N,$ also exist and, using
\eqref{S1}, we further obtain \hskip1mm
\begin{equation}\label{S2}
      \left.\begin{array}{l}\displaystyle
            \Big\|S(\dd f_n,g_n\chi_{[c,d]},P)-\int_{[c,d]}[\dd f_n]\,g_n\Big\|_Z
            \\[2mm]\hskip5mm=\displaystyle
            \Big\|S(\dd f_n,g_n\chi_{[c,d]},P)-\int_a^b[\dd f_n]\,(g_n\chi_{[c,d]})\Big\|_Z
            \\ \hskip5mm\displaystyle
            \le\Big\|\sum_{j=k+1}^{\ell}[f_n(\alpha_j)-f_n(\alpha_{j-1})]\,g_n(\xi_j)
            -\int_c^d[\dd f_n]\,g_n\Big\|_Z
            \\[7mm]\hskip10mm\displaystyle
            +\big\|[f_n(\alpha_{k+1})-f_n(c-)]\,g_n(c)\big\|_Z
            +\big\|[f_n(\alpha_{\ell+1})-f_n(d+)]\,g_n(d)\big\|_Z.
      \end{array}\right\}
\end{equation}\\
Thus, in view of \eqref{equivalence-1}, \eqref{equivalence-2}, and \eqref{S2}
and since $\{([\alpha_{j{-}1},\alpha_j],\xi_j)\,{:}\,j=k{+}1,\dots,\ell\}$
is a~$\wt{\delta}$-fine partition of $[c,d],$  we have
      \[
      \Big\|S(\dd f_n,g_n\chi_{[c,d]},P)-\int_a^b[\dd f_n]\,(g_n\chi_{[c,d]})\Big\|_Z
      <\Big(1+\frac{\|g_n(c)\|_Y}{1{+}\|g_n(c)\|_Y}+\frac{\|g_n(d)\|_Y}{1{+}\|g_n(d)\|_Y}\Big)
      \,\frac{\eps}3<\eps
      \]\smallskip
for every $n\in\N.$ This proves the equi-integrability of the sequence
$\{g_n\}$ with respect to $\{f_n\}$ on $[c,d]$ according to Definition \ref{H4-Def1}.

\skipaline

\noindent(ii) \ Let $\{g_n\chi_{[c,d]}\}$ be equi-integrable with respect to
$\{f_n\}$ on $\ab.$ In particular, all the integrals
$\int_a^b[\dd f_n]\,(g_n\chi_{[c,d]}),$ $n\,{\in}\N,$ \ exist, and for any
$\eps\,{>}\,0$ there is a gauge $\wt{\delta}$ on $\ab$ such that
\begin{equation}\label{equivalence-4}
      \left.\begin{array}{l}\displaystyle
            \quad\Big\|S(\dd f_n,g_n\chi_{[c,d]},\wt{P})-\int_a^b[\dd f_n]\,(g_n\chi_{[c,d]})\Big\|_Z
            \\[5mm]
            \mbox{for all $\wt{\delta}$-fine partitions $\wt{P}$ of $\ab$ and all $n\in\N.$}
            \end{array}
      \right\}
\end{equation}\\
Note that, by Proposition \ref{int.interval}~(iii), all the integrals
$\int_c^d[\dd f_n]\,g_n,$ $n\in\N,$ also exist.

Now, let $\eps\,{>}\,0$ and gauge $\wt{\delta}$ be given such that
\eqref{equivalence-4} is true. Put $\delta(t)=\wt{\delta}(t)$ for $t\in [c,d],$
and let $P=\{([\alpha_{j{-}1},\alpha_j],\xi_j)\}$ be an arbitrary $\delta$-fine
partition of $[c,d].$
Furthermore, let
$\wt{P}=\{([\wt{\alpha}_{j{-}1},\wt{\alpha}_j],\wt{\xi}_j)\}$ be a $\wt{\delta}$-fine
partition of $\ab,$ whose restriction to $[c,d]$ coincides with~$P.$ In other words,
there are $k,\ell\in\N$ such that
\begin{gather*}
      2\le k<\ell=k+\nu(P)<\nu(\wt{P})-2,
      \\
      \wt{\alpha}_{k+j}=\alpha_j \mbox{ \ for\ } j\in\{0,1,\dots,\ell-k\}
      \quad\mbox{and}\quad
      \wt{\xi}_{k+j}=\xi_j\mbox{ \ for\ } j\in\{1,\dots,\ell-k\}.
\end{gather*}
Furthermore, we may assume that
\begin{gather*}
      \wt{\alpha}_k=\wt{\xi}_k=\wt{\xi}_{k+1}=\alpha_0=\xi_1=c,\quad
      \wt{\alpha}_ell=\wt{\xi}_{\ell}=\wt{\xi}_{\ell+1}=\alpha_{\nu(P)}=d,
      \\
      c-\Delta<\wt{\alpha}_{k{-}1}<c \quad\mbox{and}\quad
      d<\wt{\alpha}_{\ell{+}1}<d+\Delta,
\end{gather*}
where $\Delta\in (0,\min\{b-d,c-a\})$ is given by \eqref{equivalence-2}.
Thus, considering Proposition \ref{int.interval}~(iii), for an arbitrary
$n\in\N,$ we can successively deduce
\begin{equation*}
      \left.
            \begin{array}{ll}\displaystyle
                  \Big\|S(\dd f_n,g_n,P)-\int_c^d[\dd f_n]\,g_n\Big\|_Z
                  \\[4mm]\hskip5mm\displaystyle
                  =\Big\|\sum_{j{=}1}^{\nu(P)}[f_n(\alpha_j)-f_n(\alpha_{j{-}1})]\,g_n(\xi_j)\int_c^d[\dd f_n]\,g_n\Big\|_Z
                  \\[6mm]\hskip5mm\displaystyle
                  =\Big\|S(\dd f_n,g_n\chi_{[c,d]},\wt{P})-[f_n(\wt{\alpha}_{\ell{+}1})-f_n(d)]\,g_n(d)-[f_n(\wt{\alpha}_{k{-}1})-f_n(c)]\,g_n(c)
                  \\[4mm]\hskip11mm\displaystyle
                  -\int_{[c,d]}[\dd f_n]\,g_n+\Delta^-f_n(c)\,g_n(c)+\Delta^+f_n(d)\,g_n(d)\Big\|_Z.
            \end{array}
      \right\}
\end{equation*}
Accordingly, from relations \eqref{equivalence-2} and \eqref{equivalence-4},
we finally obtain
\begin{equation*}
      \left.\begin{array}{ll}\displaystyle
            \Big\|S(\dd f_n,g_n,P)-\int_c^d[\dd f_n]\,g_n\Big\|_Z
            \\[4mm]\hskip5mm\displaystyle
            \le\Big\|S(\dd f_n,g_n\chi_{[c,d]},\wt{P})-\int_a^b[\dd f_n]\,(g_n\chi_{[c,d]})\Big\|_Z
            \\[5mm]\hskip11mm\displaystyle
            +\|f_n(c-)-f_n(\wt{\alpha}_{k{-}1})\|_X\,\|g_n(c)\|_Y
            +\|f_n(\wt{\alpha}_{\ell{+}1})-f_n(d+)\|_X\,\|g_n(d)\|_Y<\eps.
      \end{array}\right\}
\end{equation*}
This verifies relation \eqref{H3-1*}.
\endproof

\skipaline

Rather surprisingly, the same equivalences provided by
Lemma \ref{equi.properties-2A} also hold for other types of intervals.

\skipaline

\begin{lemma}\label{equi.properties-2}
Let $J$ be an arbitrary subinterval of $\ab,$ let the sequence $\{f_n\}\subset G(\ab;X)$
be equiregulated\emph{,} and let the sequence $\{g_n\}$ of functions mapping $\ab$ into $Y$
be pointwise  bounded on $\ab.$ Then\emph{,} $\{g_n\}$ is equi-integrable with respect to $\{f_n\}$
on $J$ if and only if condition \eqref{H3-1*} is satisfied.
\end{lemma}
\proof
Let $c=\inf J$ and $d=\sup J.$ For the case $J=[c,d],$ the proof was given by
the previous lemma.

\skiphalfaline

Let $J=(c,d).$ First, assume that \eqref{H3-1*} is true. Let $\eps>0$
be given and let $\delta$ be a gauge on $[c,d]$ such that
\begin{equation}\label{H3-2*}
      \left.\begin{array}{l}\displaystyle
            \quad\Big\|S(\dd f_n,g_n,P)-\int_c^d[\dd f_n]\,g_n\Big\|_{Z}<\frac\eps 3
            \\[5mm]
            \mbox{holds for any $\delta-$fine partition $P$ \ of \ } [c,d]
            \mbox{\ and any \ } n\in\N\,.
            \end{array}
      \right\}
\end{equation}
As the sequence $\{f_n\}$ is equiregulated, we may choose $\Delta\in\left(0,\frac{d-c}{2}\right)$
in such a way that the inequalities
\begin{equation}\label{Delta}
      \left.\begin{array}{l}\hskip10mm\displaystyle
            \Big\|f_n(t)-f_n(c+)\Big\|_X\,{<}\,\frac{\eps}{3\,(1{+}\sup\{\|g_n(c)\|_Y\,{:}\,n\in\N\})}\,
            \\[5mm]
            \hskip63mm\mbox{for all \ } t\in (c,c\,{+}\,\Delta)
            \\[1mm]\mbox{and}\quad
            \\[2mm]\hskip10mm\displaystyle
            \Big\|f_n(d-)-f_n(t)\Big\|_X\,{<}\,\frac{\eps}{3\,(1{+}\sup\{\|g_n(d)\|_Y\,{:}\,n\in\N\})}\,
            \\[5mm]
            \hskip63mm\mbox{for all \ } t\in (d\,{-}\,\Delta,d)
      \end{array}\right\}
\end{equation}\\
hold for each $n\in\N.$ Let $P=\{([\alpha_{j-1},\alpha_j],\xi_j)\}$ be an arbitrary
$\delta$-fine partition of $[c,d].$ Denote $m=\nu(P).$  In view of \cite[Lemma 6.2.11]{M3},
we may assume the following relation
\[
   c=\xi_1<\alpha_1<c+\Delta<d-\Delta<\alpha_{m-1}<d=\xi_m.
\]
Thus,
\begin{equation*}
      \left.\begin{array}{ll}\displaystyle
            S(\dd f_n,g_n\chi_{(c,d)},P)
            =\sum_{j=2}^{m-1}[f_n(\alpha_j)-f_n(\alpha_{j-1})]\,g_n(\xi_j)
            \\[4mm]\hskip7mm\displaystyle
            \quad=S(\dd f_n,g_n,P)-[f_n(\alpha_1)-f_n(c)]\,g_n(c)-[f_n(d)-f_n(\alpha_{m-1})]\,g_n(c)
            \\[6mm]\hskip7mm\displaystyle
            \quad=S(\dd f_n,g_n,P)-\Delta^+f_n(c)\,g_n(c)-\Delta^-f_n(d)\,g_n(d)
            \\[4mm]\hskip16mm\displaystyle
            -[f_n(\alpha_1)-f_n(c+)]\,g_n(c)-[f_n(d-)-f_n(\alpha_{m-1})]\,g_n(c).
      \end{array}\right\}
\end{equation*}\\
Accordingly, using \eqref{H3-2*}, \eqref{Delta}, and Proposition \ref{int.interval}~(i),
we get for any $n\in\N$
\begin{equation*}
      \left.\begin{array}{l}\displaystyle
            \Big\|S(\dd f_n,g_n\chi_{(c,d)},P)-\int_{(c,d)}[\dd f_n]\,g_n\Big\|_Z
            \le\Big\|S(\dd f_n,g_n,P)-\int_c^d[\dd f_n]\,g_n\Big\|_Z
            \\[4mm]\hskip10mm\displaystyle
            +\|f_n(\alpha_1)-f_n(c+)\|_X\,\|g_n(c)\|_Y
            +\|f_n(d-)-f_n(\alpha_{m-1})\|_X\,\|g_n(d)\|_Y<\eps.
      \end{array}\right\}
\end{equation*}\\
Hence, $\{g_n\}$ is equi-integrable with respect to $\{f_n\}$ on $(c,d).$

\skiphalfaline

Similarly, we would prove the reverse implication and the corresponding equivalences
for the remaining cases $J\,{=}\,[c,d)$ and $J\,{=}\,(c,d].$
\endproof

\skipaline

The next propositions summarize the properties of the equi-integrability on the subintervals
of $\ab.$

\skipaline

\begin{corollary}\label{equi.properties}
Let $J$ be an arbitrary subinterval of $\ab,$ let the sequence $\{f_n\}\subset G(\ab;X)$
be equiregulated\emph{,} and let the sequence $\{g_n\}$ of functions mapping $\ab$ into $Y$
be pointwise bounded on $\ab.$ Then\emph{,}
\begin{enumerate}\leftskip=-5mm	
\item[(i)]
If the sequence $\{g_n\}$ is equi-integrable with respect to $\{f_n\}$ on $\ab,$ then
it is equi-integrable with respect to $\{f_n\}$ on each subinterval $J$ of $\ab.$
\vskip2mm

\item[(ii)]
For every $c\in(a,b),$ the sequence $\{g_n\}$ is equi-integrable with respect to
$\{f_n\}$ on $\ab$ if and only if it is equi-integrable with respect to $\{f_n\}$
on both intervals $[a,c]$ and $[c,b].$
\end{enumerate}
\end{corollary}
\proof \ (i) \ By \cite[Proposition 8]{S}, the integral $\int_{c}^{d}[\dd f_n]\,g_n$ exists
for all $n\in\N.$ Moreover, by the Cauchy equi-integrability criterion \ref{H4-Cauchy},
hypothesis \eqref{H4-Cauchy-1} is satisfied. It can be shown in a~rather routine way that
then \eqref{H4-Cauchy-1} holds also on $[c,d].$ Therefore, by Theorem~\ref{H4-Cauchy},
$\{g_n\}$ is equi-integrable with respect to $\{f_n\}$ on $[c,d].$

\skiphalfaline

\noindent(ii) \
The necessity part of the assertion (ii) follows immediately from the first
assertion of this corollary.

On the other hand, let $\eta>0$ and $c\in(a,b)$ be given. Let $\delta_a$ and $\delta_b$
be the gauges satisfying \eqref{H3-1} for $\ab$ replaced by $[a,c]$ or $[c,b],$
respectively. Then, for a given $\eta>0,$ the condition \eqref{H3-1} will be satisfied
if we define
\[
  \delta(t)=\begin{cases}
                  \delta_a(t) &\mbox{if \ } t\in[a,c)\,,
                  \\
                  \min\{\delta_a(c),\delta_b(c)\} &\mbox{if \ } t=c\,,
                  \\
                  \delta_b(t) &\mbox{if \ } t\in(c,b]\,.
            \end{cases}
\]\vskip-6mm
\endproof

\skipaline

The next assertion is a direct consequence of Lemma~\ref{equi.properties-2}
and Corollary~\ref{equi.properties}.

\skipaline

\begin{corollary}\label{equi.cor}
Let the sequence $\{f_n\}\subset G(\ab;X)$ be equiregulated\emph{,} and let the sequence
$\{g_n\}$ of mappings of $\ab$ into $Y$ be pointwise bounded.
Then\emph{,} $\{g_n\}$ is equi-integrable with respect to $\{f_n\}$ on $\ab$ if and only
if the sequence $\{g_n\chi_{E}\}$ is equi-integrable with respect to $\{f_n\}$
on $\ab$ for every elementary subset $E$ in $\ab.$
\end{corollary}

\skiphalfaline

The well-known Saks-Henstock lemma (see e.g., \cite[Lemma 16]{S}) states that the Riemannian
sums not only approximate the integrals in the "gauge topology" over the entire interval,
but also over suitably chosen systems of subintervals. Next, we will show that
the equi-integrability implies a uniform Saks-Henstock property. However, first, let
us introduce the notion of a $\delta-$fine system.

\skipaline

\begin{definition}\label{partpart}
{\rm
The set (see e.g., \cite[Lemma 16]{S})
\[
    S=\left\{([\beta_j,\gamma_j],\xi_j)\,{:}\,j=1,2,\dots,m\right\}
\]
is a $\delta-$fine system in $\ab$ if
\[
   a\le\beta_1\le\xi_1\le\gamma_1\le\beta_2\le\xi_2\le\gamma_{2}
   \le\dots\le\beta_m\le\xi_m\le\gamma_m\le b
\]
and
\[
   [\beta_j,\gamma_j]\subset(\xi_j-\delta(\xi_j),\xi_j+\delta(\xi_j))
   \quad\mbox{for\ } j=1,2,\dots,m.
\]
Similarly, for divisions and tagged partitions of the interval $\ab,$
we will denote by $\nu(S)$ the number of the intervals contained in $S,$
i.e. $\nu(S)=m$ in the above situation.

Let $T\subset\ab$ and let
$S=\left\{\left([\beta_j,\gamma_j],\xi_j\right)\,{:}\,j\,{=}\,1,2,\dots,m\right\}$
be a $\delta-$fine system in $\ab.$
Then, $P$ is called a {\em $T-$tagged
$\delta-$fine system} in $\ab$ if $\xi_j\in T$ for every $j\,{=}\,1,2,\dots,m.$
}\end{definition}

\skipaline

\begin{proposition}
Let $\{f_n\}\subset G(\ab,X)$ be equiregulated\emph{,} and let the sequence $\{g_n\}$
of functions mapping $\ab$ into $Y$ be equi-integrable with respect to $\{f_n\}$
on $\ab.$ Furthermore\emph{,} let $\eps>0$ be given arbitrarily and let $\delta$ be a gauge on $\ab$ such that
\begin{align*}
      \Big\|S(\dd f_n,g_n,P)-\int_a^b[\dd f_n]\,g_n\Big\|_{Z}<\eps
\end{align*}
for every $\delta-$fine partition $P$ of $\ab$ and every  $n\in\N.$ Then\emph{,}
\begin{align*}
      \Big\|\sum_{j{=}1}^{\nu(S)}\Big([f_n(\gamma_j)-f_n(\beta_j)]\,g_n(\xi_j)
      -\int_{\beta_j}^{\gamma_j}[\dd f_n]\,g_n\Big)\Big\|_{Z}\le\eps
\end{align*}
holds for every $\delta-$fine system $S=\{([\beta_j,\gamma_j],\xi_j)\}$ in $\ab$
and every $n\in\N.$
\end{proposition}
\proof
Assume that the system $\{([\gamma_j,\beta_j],\xi_j):j=1,2,{\dots},n\}$
satisfies the following conditions
\begin{gather*}\label{sys}
      \left.\begin{array}{l}
            a\le \gamma_1\le\xi_1\le\beta_1\le\gamma_2\le\dots\le\gamma_n\le\xi_n\le\beta_n\le b,
            \\[2mm]
            [\gamma_j,\beta_j]\subset(\xi_j-\delta(\xi_j),\xi_j+\delta(\xi_j))
        	\quad\mbox{for\ } j=1,\dots,n,
      \end{array}\right\}
\end{gather*}
and denote $\gamma_0=a$ and $\beta_{n+1}=b.$ Then, by Corollary \ref{equi.cor},
the sequence $\{g_n\}$ is equi-integrable with respect to $\{f_n\}$ on each
subinterval $[\beta_j,\gamma_{j+1}].$ Hence, we can apply the method of the proof
of the Saks-Henstock lemma in \cite[Lemma 6.5.1]{M3} to complete the proof of
this proposition.
\endproof

\skipaline\vskip-0.5mm

\section{Harnack extension principle and its applications}\label{HarnackApp}

\skiphalfaline

Assume that
\begin{equation}\label{intro_Harnack_1}
      \left.\begin{array}{l}\displaystyle
            (a,b)\setminus T=\bigcup_{i\in\N}(a_i,b_i),\quad\mbox{where\ } (a_i,b_i), i\in\N,
            \\[5mm]
            \mbox{are pairwise disjoint open subintervals of} \,(a,b).
      \end{array}\right\}
\end{equation}
For every $i\in\N,$ let $J_i$ denote $\left[a_i,b_i\right],$ $\left[a_i,b_i\right),$ $\left(a_i,b_i\right],$ or $\left(a_i,b_i\right);$ let $E$ be an ordered finite union of $J_i,$ i.e., $E=\bigcup_{i=p}^{p+q} J_i$ for some $p,q\in\N,$
and let the (Kurzweil-Henstock) integrals $\int_{J_i} g\,\dd t$ or $\int_{a_i}^{b_i} g\,\dd t$ exist for each $i\in\N.$
The following equation
\begin{equation}\label{intro_Harnack_2}
      \int_{\ab} g\,\dd t=\int_a^b g\,\dd t = \int_{(a,b)} g\,\dd t
\end{equation}
(provided with one of the integrals $\int_{\ab} g\,\dd t,$ $\int_a^b g\,\dd t,$ or $\int_{(a,b)} g\,\dd t$ exists)
and the situations of \eqref{differ-equal-2} and \eqref{rem.inf.sup}, which are true for the Kurzweil-Henstock
integral, have verified the statements below
\begin{equation}\label{intro_Harnack_3}
      \int_{\ab\setminus T} g\,\dd t=\int_{(a,b)\setminus T} g\,\dd t \, \, \, \Big(\text{in case} \, \int_{T} g\,\dd t \, \text{exists}\Big),
\end{equation}
\begin{equation}\label{intro_Harnack_4}
      \int_{J_i} g\,\dd t=\int_{a_i}^{b_i} g\,\dd t, \, \, \, \text{for each} \, i\in\N,
\end{equation}
and
\begin{equation}\label{intro_Harnack_5}
      \int_{E} g\,\dd t=\int_{\bigcup_{i=p}^{p+q} J_i} g\,\dd t=\sum_{i=p}^{p+q}\int_{a_i}^{b_i} g\,\dd t,
\end{equation}
respectively. As a consequence, under the conditions in Theorem \ref{harnack}, for the Kurzweil-Henstock
integral we have that the integrals
\[
      \int_a^b g\,\dd t,\,  \int_{(a,b)\setminus T} g\,\dd t,  \, \, \, \text{and} \, \, \, \int_{\bigcup_{i\in\N} J_i} g\,\dd t,
\]
exist and the following relation
\[
      \int_a^b g\,\dd t=\int_a^b g\chi_{T}\,\dd t\,{+}\int_{\bigcup_{i\in\N}J_i}\hskip-4mm g\,\dd t
      =\int_a^b g\chi_{T}\,\dd t\,{+}\sum_{i=1}^{\infty}\int_{a_i}^{b_i}g\,\dd t
\]
holds.

\skipaline

The Harnack extension is also valid for the Kurzweil-Henstock integrable real-valued functions defined
on measure spaces endowed with locally compact metric topologies, (see e.g., \cite[Theorem 5.1]{NWL}).
It applies two important concepts, namely, the $\delta-$fine cover and nonabsolute subset, involving
the integration over elementary sets and fulfilling \eqref{intro_Harnack_2}, \eqref{intro_Harnack_3},
\eqref{intro_Harnack_4}, and \eqref{intro_Harnack_5}.

\skipaline

Nonetheless, in view of \eqref{int.closed.interval-added}, \eqref{r3,1},
\eqref{differ-equal-1}, or Remark \ref{min.decomp} (ii), the relations \eqref{intro_Harnack_2},
\eqref{intro_Harnack_3}, \eqref{intro_Harnack_4}, and \eqref{intro_Harnack_5} need not be true
for the KS-integral. Hence, neither Theorem \ref{harnack} nor \cite[Theorem 5.1]{NWL} could be simply
extended to the KS-integration.

\skipaline

We need the following definition to construct the Harnack extension principle
for the KS-integral.

\skiphalfaline

\begin{definition}\label{pre_Harnack3}
{\rm\
Let $E$ be an elementary set in $\R$ and let $T$ be a closed subset of $\overline{E}.$
The sequence $\{E_i\}$ of elementary subsets of $E\setminus T$ is called a \emph{proper cover}
of $E\setminus T$ if $E\setminus T=\bigcup_{i\in\N} E_i.$
}\end{definition}

\skiphalfaline

\begin{remark}\label{pre_Harnack4}
{\rm Let $E$ be an elementary set in $\R$ and let $T$ be a closed subset of $E.$
\begin{itemize}
\item[(i)]
Then by \cite[Lemma 5.1]{NWL}, there must exist elementary sets $E_1, E_2,\dots$ such that
$E_i\subseteq E\setminus T$ for every $i\in\N$ and
      \begin{equation}\label{pre_Harnack2}
      \mu\Big(\left(E\setminus T\right)\setminus \bigcup_{i\in\N} E_i\Big)=0
      \end{equation}
($\mu$ is a signed measure). The sequence $\left\{E_i\right\}$ fulfilling \eqref{pre_Harnack2} in such a way that
$\left(E\setminus T\right)\setminus \bigcup_{i\in\N} E_i=\emptyset,$ is a proper cover
of $E\setminus T.$

\item[(ii)]
Consider $E=(a,b)\subset\ab$ and a closed subset $T\subset\ab.$ Then, $(a,b)\setminus T$
is an open set and thus can be written as the union of a countable of pairwise disjoint
intervals $(a_i,b_i)\subset (a,b)\setminus T,$ $i\in\N.$ Hence, the sequence
$\left\{(a_i,b_i)\right\}$ is a proper cover of $(a,b)\setminus T$ (see e.g., Theorem \ref{harnack}
or \eqref{intro_Harnack_1}).
\end{itemize}
}
\end{remark}

\skipaline

Applying Theorem \ref{H3-convergence} and Definition \ref{pre_Harnack3} together, we are now ready to prove the Harnack extension principle for the KS-integral by taking $E_i,\,i\in\N$ to be mutually disjoint elementary sets in $\ab,$ such that $\ab\setminus T=\bigcup_{i\in\N}E_i$ for a closed subset $T$ of $\ab.$

\skipaline

\begin{theorem}\label{Harnack}{\rm({Harnack extension principle})}
Let $f\in(\mathcal{B})G(\ab,X)$ and $g:\ab\to Y.$ Let $T$
be a closed subset of $\ab$ and
\[
   S:=\ab\setminus T=\bigcup_{i\in\N}E_i,
\]
where $E_i,$ $i\in\N,$ are mutually disjoint elementary sets in $\ab.$ Furthermore\emph{,} put
$S_n=\bigcup_{i=1}^n E_i$ for $n\in\N$ and assume that the integral $\int_{T}\f\,g$ exists
and the sequence $\{g\chi_{S_n}\}$ is equi-integrable with respect to $f$ on $\ab.$
Then the integrals $\int_{S}\f\,g,$ $\int_{\ab}\f\,g,$ and
$\int_{E_i}\f\,g$ exist for all $i\in\N.$ Moreover\emph{,}
\begin{align}\label{H1}
      \int_{\ab}\f\,g&=\int_T\f\,g+\int_S\f\,g,
      \\\noalign{\noindent\mbox{where}}\label{H2}
      \int_S\f\,g&=\sum_{i=1}^{\infty}\int_{E_i}\f\,g.
\end{align}
\end{theorem}
\proof
Clearly, $S_n$ is an elementary set in $\ab$ for every $n\,{\in}\,\N$ and $E_i\subset S_n$ for every $i=1,2,\,\dots\,n.$
Since the integral $\int_{S_n}\f \,g=\int_a^b\f\,(g\chi_{S_n})$ exists for every $n\,{\in}\,\N,$ $E_i$ for $i=1,2,\,\dots\,n$ are mutually disjoint, and $f\in(\mathcal{B})G(\ab,X),$ by Theorem \ref{int.elem.set}  $\text{(i)},$ Corollary \ref{C1}, and Proposition \ref{P2berubah}, the integrals $\int_{E_i}\f \,g,$ for $i=1,2,\dots,n,$ exist and
\begin{align}\label{int.series}
      \int_{S_n}\f\,g=\sum_{i=1}^{n}\int_{E_i}\f\,g
      \quad\mbox{for any \ } n\in\N.
\end{align}
Furthermore, as
\begin{align}\label{int.series.2}
      (g\chi_{S})(t)=\lim_{n\to\infty}(g\chi_{S_n})(t)
      \quad\mbox{for all \ } t\in\ab,
\end{align}
making use of Theorem \ref{H3-convergence} and \eqref{int.series}, we obtain
\[
     \int_S\f\,g=\lim_{n\to\infty}\int_{S_n}\f\,g
                =\sum_{i=1}^{\infty}\int_{E_i}\f\,g,
\]
i.e., \eqref{H2} is true. Finally, \eqref{H2} together with Proposition \ref{P1} and Corollary \ref{C1}
imply \eqref{H1}. This completes the proof.
\endproof

\skipaline

\begin{remark}{\rm
Recall, cf. Remark \ref{R}, that under the assumptions of Theorem \ref{Harnack},
the integrals $\int_{\ab}\f\,g$ and $\int_a^b\f\,g$  coincide.
}\end{remark}

\skipaline

If $T\subset\ab,$ then, in general, the existence of the integral
$\int_{\ab}\f \,g$ does not imply the existence of the integral $\int_{T}\f\,g$
(except when $T$ is an elementary subset of $\ab$; see Proposition \ref{int.elem.set} (i)).
This is demonstrated by the following example. Let $g:[0,1]\to\R$ be defined by
\[
   g(t)=\begin{cases}
              0 &\mbox{if \ } t=0,
              \\[1mm]
              2t\,\cos\frac{\pi}{t^2}+\frac{2\pi}{t}\,\sin\frac{\pi}{t^2}
              &\mbox{if \ } 0<t\le 1.
        \end{cases}
\]
Then, $g$ is Kurzweil-Henstock integrable on $[0,1].$ However, the integral
$\int_{T}g\,\dd t=\infty$ if $T\,{:=}\,\{t\in [0,1]: g(t)\ge 0\}$ (see e.g., \cite{K2},
\cite{L1}). This means that the opposite implication given by the previous
theorem does not hold, in general. The next theorem provides a certain affirmative
result for the case that the integrator $f$ is simply-regulated.

\skiphalfaline

\begin{theorem}\label{int.any.subset}
Let $f\in(\mathcal{B})G(\ab,X)$ and $g:\ab\to Y.$ Let $T$ be a closed subset of $\ab$ and
\[
   S:=\ab\setminus T=\bigcup_{i\in\N}E_i,
\]
where $E_i,$ $i\in\N,$ are mutually disjoint elementary sets in $\ab.$ Furthermore\emph{,} put
$S_n=\bigcup_{i=1}^n E_i$ for $n\in\N$ and assume that the integral $\int_{\ab}\f\,g$ exists
and the sequence $\{g\chi_{S_n}\}$ is equi-integrable with respect to $f$ on $\ab.$ Then\emph{,}
\begin{itemize}
      \item[{\rm(i)}]
           All the integrals
           $\int_{E_i}\f\,g$ for $i\in\N,$ $\int_{S}\f\,g,$ and $\int_{T}\f\,g$ exist
           and the equalities
                \[
                \int_{S}\f\,g=\sum_{i=1}^{\infty}\int_{E_i}\f\,g
                \mbox{ \ and \ }
                \int_{\ab}\f\,g=\int_{T}\f\,g+\int_{S}\f\,g
                \]
           are true.\vskip2mm
      \item[{\rm(ii)}]
           For every $\eps>0$, there is a gauge $\delta$ on $T,$ such that
                \[
                \Big\|\sum_{j{=}1}^{\nu(Q)}\Big(\int_{\alpha_j}^{\beta_j}\f\,g\chi_{T}
                -\int_{\alpha_j}^{\beta_j}\f\,g\Big)\Big\|_Z<\eps
                \]
           for every $T-$tagged $\delta-$fine system $Q=\{([\alpha_j,\beta_j],\xi_j)\}$ in $\ab.$
\end{itemize}
\end{theorem}
\proof \ (i)\quad
As a result of the equi-integrability property of the sequence $\{g\chi_{S_n}\}$ with
respect to $f\in(\mathcal{B})G(\ab,X),$ as in the proof of Theorem \ref{Harnack}, we can show that all
the integrals $\int_{E_i}\f\,g,$ $i\in\N,$ exist and
\[
  \int_{S}\f\,g=\sum_{i=1}^{\infty}\int_{E_i}\f\,g.
\]
Moreover, $S\cup T=\ab$ and $S\cap T=\emptyset.$ Therefore, by
Proposition \ref{P1} and Corollary \ref{C1}, where we put $S_1=S$ and $S_2=T,$ we also obtain that
the integral $\int_{T}\f\,g$ exists and the equality
\[
   \int_{T}\f\,g =\int_{\ab}\f\,g-\int_{S}\f\,g
\]
holds.

\noindent(ii)\quad
Let $\eps>0$ be given and let $\delta$ be a gauge on $\ab$ such that
\begin{gather*}
      \Big\|S(\dd f,g,P)-\int_a^b\f\,g\Big\|_Z<\frac{\eps}{6}
      \\\noalign{\noindent\mbox{and}}
      \Big\|S(\dd f,g\chi_{T},P)-\int_a^b\f \,(g\chi_{T})\Big\|_Z<\frac{\eps}{6}
\end{gather*}
whenever $P$ is a $\delta-$fine partition of $\ab.$

Suppose now that $Q=\{([\alpha_j,\beta_j],\xi_j)\}$ is a $T-$tagged $\delta-$fine
system in $\ab.$ Then, using the Saks-Henstock lemma (see \cite[Lemma 16]{S}), we
deduce
\begin{equation*}
      \left.\begin{array}{ll}\displaystyle
            \Big\|\sum_{j=1}^{\nu(Q)}\Big(\int_{\alpha_j}^{\beta_j}\f\,(g\chi_{T})
                                -\int_{\alpha_j}^{\beta_j}\f\,g\Big)\Big\|_Z
            \\[4mm]\hskip6mm\displaystyle
            \le\Big\|\sum_{j=1}^{\nu(Q)}\Big(\int_{\alpha_j}^{\beta_j}\f\,(g\chi_{T})
                               -[f(\beta_j)-f(\alpha_j)]\,(g\chi_{T})(\xi_j)\Big)\Big\|_Z
            \\[5mm]\hskip13mm\displaystyle
            +\Big\|\sum_{j{=}1}^{\nu(Q)}\Big([f(\beta_j)-f(\alpha_j)]\,g(\xi_j)
                                  -\int_{\alpha_j}^{\beta_j}\f\,g\Big)\Big\|_Z<\eps.
      \end{array}\right\}
\end{equation*}
\noindent
This completes the proof.
\endproof

\skipaline

The next result is essentially a corollary of Theorem \ref{int.any.subset}.

\skipaline

\begin{theorem}\label{the-last}
Let $f\in G(\ab,X)$ and $g:\ab\to Y.$ Let $E$ be an elementary set in $\ab,$ $T$ be a closed subset of $E,$ and
\[
   \ab\setminus T=\bigcup_{i\in\N}E_i,
\]
where $E_i,$ $i\in\N,$ are mutually disjoint elementary sets in $\ab.$
Furthermore\emph{,} put
$S_n=\bigcup_{i=1}^n E_i$ for $n\in\N$ and assume that the integral $\int_{E}\f\,g$ exists
and the sequence $\{g\chi_{S_n}\}$ is equi-integrable with respect to $f$ on $\ab.$
Then\emph{,} the integral $\int_{T}\f \,g$ exists.
\end{theorem}
\begin{proof}
Without loss of generality, by Remark \ref{min.decomp}, throughout the proof, we may assume that $E=\bigcup_{k=1}^m J_k$ where $\{J_1, J_2,\dots,J_m\}$ is a minimal decomposition of $E.$
Hence, by the hypothesis, the integrals $\int_{J_k}\f \,g,$ for $k=1,2,\dots,m,$ exist.
For every $k=1,2,\dots,m,$ we denote $c_k=\inf(J_k)$ and $d_k=\sup(J_k)$ and set
      \begin{gather*}
      T_k=T\cap [c_k,d_k].
      \end{gather*}
Evidently, $\bigcup_{k=1}^m T_k=T$
and $\bigcap_{k=1}^m T_k=\emptyset$ due to the minimal decomposition of $\{J_1, J_2,\dots,J_m\}.$
Then, from Remark \ref{rem3.5}, we obtain that for every $k=1,2,\dots,m,$ the integrals
      \begin{align*}
      \int_{(c_k,d_k)}\f\,g,\,\,\int_{[c_k,d_k)}\f\,g,\,\,\int_{(c_k,d_k]}\f\,g,\,\,
      \int_{[c_k,d_k]}\f\,g,\,\, \text{and} \,\, \int_{c_k}^{d_k}\f\,g
      \end{align*}\noindent
exist.

Furthermore, for every subset $[c_k,d_k]\subset \ab,$ $k=1,2,\dots,m,$ we may construct a sequence
                       \[
                       \left\{F_{i}^{(k)}=E_i\cap [c_k,d_k];\,i\in\N\right\}
                       \]
                   of mutually disjoint elementary sets in $[c_k,d_k],$ such that
                       \[
                       \bigcup_{i\in\N}F_{i}^{(k)}=[c_k,d_k]\setminus T_k=[c_k,d_k]\setminus T,\,\,\text{for}\,k=1,2,\dots,m.
                       \]
                   Moreover, if we put $H_{n}^{(k)}=\bigcup_{i=1}^n F_{i}^{(k)}$ for $n\in\N,$ then we will have $H_{n}^{(k)}\subseteq S_n$ for every $n\in\N.$

                   From \eqref{int.series.2}, in the proof of Theorem \ref{Harnack}, the sequence $\{g\chi_{S_n}\}$ converges
                   pointwise on $\ab,$ which further implies the boundedness of the sequence
                   $\{(g\chi_{S_n})(t)\}$ in $Y$ for each $t\in\ab.$ Accordingly, we obtain that the sequence
                   $\{(g\chi_{H_{n}^{(k)}})(t)\}$ is bounded in $Y$ for each $t\in\ab$ and $k=1,2,\dots,m.$
                   As a consequence, by Proposition \ref{equi.properties} (i) and Corollary \ref{equi.cor},
                   the sequence
                   $\{g\chi_{H_{n}^{(k)}}\}$ is equi-integrable with respect to $f$ on $[c_k,d_k],$ for every $k=1,2,\dots,m.$
                   Therefore, by Theorem \ref{int.any.subset}, the integrals
                   $\int_{\bigcup_{i\in\N}F_{i}^{(k)}}\f\,g$ and $\int_{T_k}\f\,g$ exist
                   and
                           \[
                           \int_{T_k}\f\,g=\int_{[c_k,d_k]}\f\,g-\int_{\bigcup_{i\in\N}F_{i}^{(k)}}\f\,g
                           \]
                   holds for every $k=1,2,\dots,m.$ Finally, the existence of the integral $\int_{T}\f\,g$ is obtained by making use of Remark \ref{R3.8} and Proposition \ref{P2berubah}.
\end{proof}

\skipaline

\section{Conclusion}
In the present study, using known results for the Kurzweil-Stieltjes integrals, we were able to prove theorems concerning the equi-integrability and equiregulatedness involving elementary sets, to construct the Harnack extension principle in the setting of the Kurzweil-Stieltjes integration, and to show its essential role in solving a problem on how to get the existence of integral $\int_{T}\f\,g$ for an arbitrary closed subset $T$ of the elementary set $E$ if the integral $\int_{E}\f\,g$ exists.

\skipaline

\noindent \textbf{Acknowledgments} The author is very grateful to Professor Milan Tvrd\'y (Czech Academy of Sciences, Czech Republic) for his comments, suggestions, and help to finish this paper. 

\skipaline

\bibliographystyle{amsplain}

\end{document}